
\def\LaTeX{%
  \let\Begin\begin
  \let\End\end
  \let\salta\relax
  \let\finqui\relax
  \let\futuro\relax}

\def\UK{\def\our{our}\let\sz s}
\def\USA{\def\our{or}\let\sz z}

\UK



\LaTeX

\USA


\salta

\documentclass[twoside,12pt]{article}
\setlength{\textheight}{24cm}
\setlength{\textwidth}{16cm}
\setlength{\oddsidemargin}{2mm}
\setlength{\evensidemargin}{2mm}
\setlength{\topmargin}{-15mm}
\parskip2mm


\usepackage[usenames,dvipsnames]{color}
\usepackage{amsmath}
\usepackage{amsthm}
\usepackage{amssymb}
\usepackage[mathcal]{euscript}

\usepackage[T1]{fontenc}
\usepackage[latin1]{inputenc}
\usepackage[english]{babel}
\usepackage[babel]{csquotes}

\usepackage{latexsym}
\usepackage{graphicx}
\usepackage{mathrsfs}
\usepackage{mathrsfs}
\usepackage{hyperref}
\usepackage{pgfplots}
%
%
\usepackage{cite}
%
%
%


\definecolor{viola}{rgb}{0.3,0,0.7}
\definecolor{ciclamino}{rgb}{0.5,0,0.5}
\definecolor{rosso}{rgb}{0.8,0,0}
\definecolor{darkgreen}{rgb}{0,0.5,0}

\def\pier #1{{\color{blue}#1}}
\def\juerg #1{{\color{red}#1}}
\def\lastrev #1{{\color{rosso}#1}}
\def\andrea #1{{\color{red}#1}}
\def\andrealast #1{{\color{red}#1}}
\def\last #1{{\color{rosso}#1}}
\def\lastjuerg #1{{\color{blue}#1}}

\def\pier #1{#1}
\def\juerg #1{#1}
\def\lastrev #1{{#1}}
\def\andrea #1{#1}
\def\andrealast #1{#1}
\def\last #1{{#1}}
\def\lastjuerg #1{#1}



\bibliographystyle{plain}


%
\newtheorem{theorem}{Theorem}[section]

\finqui

\def\Bthm{\Begin{theorem}}
\def\Ethm{\End{theorem}}

\def\Bcor{\Begin{corollary}}
\def\Ecor{\End{corollary}}
\def\Brem{\Begin{remark}\rm}
\def\Erem{\End{remark}}

\def\Bdim{\Begin{proof}}
\def\Edim{\End{proof}}
\def\Bcenter{\Begin{center}}
\def\Ecenter{\End{center}}
\let\non\nonumber




\def\step #1 \par{\medskip\noindent{\bf #1.}\quad}


\def\Lip{Lip\-schitz}
\def\Holder{H\"older}
\def\Frechet{Fr\'echet}

\def\wk{well-known}
\def\socal{so-called}
\def\lhs{left-hand side}
\def\rhs{right-hand side}
\def\sfw{straightforward}

\def\CH{Cahn--Hilliard}



\def\multibold #1{\def\arg{#1}%
  \ifx\arg\pto \let\next\relax
  \else
  \def\next{\expandafter
    \def\csname #1#1#1\endcsname{{\bf #1}}%
    \multibold}%
  \fi \next}

\def\pto{.}

\def\multical #1{\def\arg{#1}%
  \ifx\arg\pto \let\next\relax
  \else
  \def\next{\expandafter
    \def\csname cal#1\endcsname{{\cal #1}}%
    \multical}%
  \fi \next}


\def\multimathop #1 {\def\arg{#1}%
  \ifx\arg\pto \let\next\relax
  \else
  \def\next{\expandafter
    \def\csname #1\endcsname{\mathop{\rm #1}\nolimits}%
    \multimathop}%
  \fi \next}

\multibold
qwertyuiopasdfghjklzxcvbnmQWERTYUIOPASDFGHJKLZXCVBNM.

\multical
QWERTYUIOPASDFGHJKLZXCVBNM.

\multimathop
ad dist div dom meas sign supp .


\def\accorpa #1#2{\eqref{#1}--\eqref{#2}}
\def\Accorpa #1#2 #3 {\gdef #1{\eqref{#2}--\eqref{#3}}%
  \wlog{}\wlog{\string #1 -> #2 - #3}\wlog{}}

\def\Accorpadue #1#2 #3 {\gdef #1{\ref{#2}--\ref{#3}}%
  \wlog{}\wlog{\string #1 -> #2 - #3}\wlog{}}


\def\infess{\mathop{\rm inf\,ess}}
\def\supess{\mathop{\rm sup\,ess}}

\def\graffe #1{\mathopen\{#1\mathclose\}}

\def\<#1>{\mathopen\langle #1\mathclose\rangle}
\def\norma #1{\mathopen \| #1\mathclose \|}

\def\iot {\int_0^t}
\def\ioT {\int_0^T}
\def\intQt{\int_{Q_t}}
\def\intQ{\int_Q}
\def\iO{\int_\Omega}

\def\intQtT{\int_{Q_t^T}}

\def\dt{\partial_t}
\def\dn{\partial_n}

\def\checkmmode #1{\relax\ifmmode\hbox{#1}\else{#1}\fi}

\def\aeQ{\checkmmode{a.e.\ in~$Q$}}

\def\aaQ{\checkmmode{for a.a.~$(x,t)\in Q$}}


\def\erre{{\mathbb{R}}}

\def\enne{{\mathbb{N}}}




\def\genspazio #1#2#3#4#5{#1^{#2}(#5,#4;#3)}
\def\spazio #1#2#3{\genspazio {#1}{#2}{#3}T0}

\def\L {\spazio L}
\def\H {\spazio H}

\def\C #1#2{C^{#1}([0,T];#2)}


\def\Lx #1{L^{#1}(\Omega)}
\def\Hx #1{H^{#1}(\Omega)}

\def\Ldue{\Lx 2}

\def\Huno{\Hx 1}
\def\Hdue{\Hx 2}



\let\theta\vartheta
\let\eps\varepsilon
\let\phi\varphi
\let\lam\lambda

\let\TeXchi\chi                         
\newbox\chibox
\setbox0 \hbox{\mathsurround0pt $\TeXchi$}
\setbox\chibox \hbox{\raise\dp0 \box 0 }
\def\chi{\copy\chibox}



%
%

\def\bQ{\gamma_2}
\def\bS{\gamma_4}
\def\bO{\gamma_1}
\def\bG{\gamma_3}
\def\bu{\gamma_5}
\def\bw{\gamma_6}

\def\rmin{r_-}
\def\rmax{r_+}

\def\phQ{\phi_Q}

\def\phO{\phi_\Omega}

\def\sQ{\s_Q}
\def\sO{\s_\Omega}

\def\Uad{\calU_{\ad}}
\def\uopt{\overline u}
\def\wopt{\overline w}

\def\Vp{V^*}

\def\normaV #1{\norma{#1}_V}
\def\normaH #1{\norma{#1}_H}

\let\hat\widehat

\def\Pi{\hat\pi}



\def\s{\sigma}  
\def\ss{\sigma_s}  
\def\m{\mu}	    
\def\ph{\phi}	
\def\a{\alpha}	
\def\b{\beta}	
\def\d{\delta}  
\def\et{\eta}   
\def\th{\theta} 
\def\bph{\overline\ph}  
\def\bm{\overline\m}    
\def\bs{\overline\s}    
\def\hph{\hat\ph}  
\def\hm{\hat\m}    
\def\hs{\hat\s}    
\def\z{\zeta}      
\def\ch{\chi}      
\def\ps{\psi}      
\def\J{{\cal J}}		
\def\Jred{{\cal J}_{\rm red}} 
\def\S{{\cal S}}   
\def\I2 #1{\int_{Q_t}|{#1}|^2}
\def\IN2 #1{\int_{Q_t}|\nabla{#1}|^2}
\def\IO2 #1{\iO |{#1(t)}|^2}
\def\INO2 #1{\iO |\nabla{#1}(t)|^2}
\def\UR{{\cal U}_R}
\def\CP{\boldsymbol{(CP)}}
\def\hinf{h_\infty}


\Begin{document}



\title{Optimal control of a phase field system\\[0.3cm]
	modelling tumor growth with \andrea{chemotaxis}\\[0.3cm]
	\andrea{and} singular potentials}
\author{}
\date{}
\maketitle

\Bcenter
\vskip-1cm
{\large\sc Pierluigi Colli$^{(1)}$}\\
{\normalsize e-mail: {\tt pierluigi.colli@unipv.it}}\\[.25cm]
{\large\sc Andrea Signori$^{(2)}$}\\
{\normalsize e-mail: {\tt andrea.signori02@universitadipavia.it}}\\[.25cm]
{\large\sc J\"urgen Sprekels$^{(3)}$}\\
{\normalsize e-mail: {\tt sprekels@wias-berlin.de}}\\[.5cm]
$^{(1)}$
{\small Dipartimento di Matematica ``F. Casorati'', Universit\`a di Pavia}\\
{\small via Ferrata 5, 27100 Pavia, Italy}\\[.2cm]
$^{(2)}$
{\small Dipartimento di Matematica e Applicazioni, Universit\`a di Milano--Bicocca}\\
{\small via Cozzi 55, 20125 Milano, Italy}\\[.2cm]
$^{(3)}$
{\small Department of Mathematics}\\
{\small Humboldt-Univesit\"at zu Berlin }\\
{\small Unter den Linden 6, 10099 Berlin, Germany}\\
{\small and}\\
{\small Weierstrass Institute for Applied Analysis and Stochastics}\\
{\small Mohrenstrasse 39, 10117 Berlin, Germany}

\Ecenter

%
%

%
\Begin{abstract}\noindent
A distributed optimal control problem for an extended model of phase field
type for tumor growth is addressed.
In this model, the chemotaxis effects are also taken into account.
The \juerg{control} is realized by two control variables
that design the dispensation of some drugs to the patient.
The cost functional is \juerg{of} tracking type, whereas
the potential setting has been kept quite general in order to allow 
regular and singular potentials to be considered. 
In this direction, some relaxation terms have been introduced
in the system. We show the well-posedness of the state system,
the \Frechet\ differentiability of the control-to-state operator 
in a suitable functional \juerg{analytic} framework, and, lastly, we
characterize the first-order necessary \juerg{conditions
of} optimality in terms of a variational inequality
involving the adjoint variables.

\End{abstract}

\pagestyle{myheadings}
\newcommand\testopari{\sc Colli -- Signori -- Sprekels}
\newcommand\testodispari{\sc \pier{Optimal} control for a tumor growth model}
\markboth{\testopari}{\testodispari}

\salta
\finqui

\newpage
\section{Introduction}
\label{SEC_INTRODUCTION}
\setcounter{equation}{0}
After realizing that \juerg{tumor
cells, like any other} material, have to obey physical laws,
a significant number of models have been introduced, since from a modelling
viewpoint a tumor mass \juerg{does not behave that different from other} special materials
investigated by scientists 
(see \cite{CL} and also \cite{CLLW,OHP,HDPZO,HKNZ,WLFC,WZZ}).
As \juerg{far as} diffuse interface models are concerned, we can identify two main classes. 
The first one considers the tumor and healthy cells as 
\juerg{inertialess fluids and 
includes effects generated by the fluid flow development
by} postulating a Darcy law or a Stokes--Brinkman law.
\juerg{In this connection}, we refer to
\cite{WLFC,GLSS,DFRGM ,GARL_1,GARL_4,GAR, EGAR, FLRS,GARL_2, GARL_3}, 
where \juerg{also} further mechanisms such as chemotaxis and active transport 
are taken into account.
The other class, to which our model belongs, neglects the velocity.

In this framework, let us take $\Omega\subset\erre^3$ as an open, bounded, and connected set
with smooth boundary $\Gamma$; \juerg{moreover, we set, for $0<t<T$, 
$$Q:=\Omega\times (0,T),\quad \Sigma:=\Gamma\times (0,T),\quad Q_t:=\Omega\times (0,t),
\quad Q_t^T:=\Omega\times (t,T).
$$
}
 The initial-boundary problem \juerg{under investigation then reads} as follows.
\begin{align}
	\label{State_first}
	&\a \dt \m
	+ \dt \phi- \Delta \m
	\,=\,
	(P\s -A - u) h(\ph)
	 &\hbox{in $\, Q,$}
	\\[0.5mm]
	\label{State_second}
	&\m
	\,=\,
	\b\dt \ph
	-\Delta \ph
	+ F'(\ph)
	- \ch \s
 	 &\hbox{in $\,Q,$}
	\\[0.5mm]
	\label{State_third}
	&\dt \s
	- \Delta \s
	\,=\,
	\andrea{- \ch \Delta \ph +}
	B(\ss - \s)
	- D \s h(\ph)
	+w
	 &\hbox{in $\,Q,$}
	\\[0.5mm]
	&\dn \m\,=\,\dn \ph\,=\,\dn \s \,=\,0
 	 &\hbox{on $\,\Sigma,$}
	\label{State_BC}
	\\[0.5mm] 
	&\m(0)\,=\, \m_0, \ \ph(0)\,=\,\ph_0, \ \s(0)\,=\,\s_0
  	 &\hbox{in $\,\Omega$,}
	\label{State_IC}
\end{align}
\Accorpa\Stsys State_first State_IC
where the symbol $\dn$ \juerg{indicates} the outward normal derivative \juerg{to $\Gamma$}.
The above state system consists of an extended
Cahn--Hilliard type system for the tumor phase
coupled with a reaction-diffusion equation for an unknown species acting as a nutrient \last{(e.g., glucose, oxygen, carbohydrates)}.
The system \Stsys\ is a simplification and relaxed version of the model originally proposed 
in \cite{GLSS}. Indeed,
the velocity contributions 
are neglected, and two relaxation terms are added. This choice
will allow us to consider more general potentials
that may exhibit \juerg{a} singular behavior.
By assuming different linear phenomenological laws for chemical reactions,
a different thermodynamically consistent model was introduced in
\cite{HDZO} (see also \cite{HKNZ,CLLW,HDPZO,OHP}), and the corresponding
mathematical investigations have been addressed
in \cite{FGR,CGH,CGRS_VAN,CGRS_ASY}.
In \cite{CGH,CGRS_VAN,CGRS_ASY} the same two relaxation
terms $\a\dt\m$ and $\b\dt\ph$ have been introduced.
As in the current case, their presence 
allowed the authors to take into account more general potentials
that may be singular and also nonregular.
Moreover, in \cite{CGRS_VAN,CGRS_ASY}, the authors
\juerg{pointed} out how $\a$ and $\b$ can be \juerg{set to zero, by} providing
the proper framework in which a limit system can be identified and uniquely solved
\last{(let us also mention 
\cite{CGS_frac_asy, CGS_frac_wp, CGS_frac_control}, where similar problems have been addressed for the
case of fractional operators)}. 
Next, we mention \cite{FRL}, where a similar \juerg{nonlocal} version was studied
for the case of singular potentials and degenerate mobilities.
Let us also point out \cite{MRS,CRW}, where the long-time behavior of these 
models \juerg{was} analyzed in terms of the convergence to
equilibrium and \juerg{of} the existence of a global attractor, respectively.
For further models, discussing the case of multispecies,
we address the reader to \cite{DFRGM,FLRS}.

Now, let us briefly describe the role of the occurring
symbols from a \juerg{modeling} viewpoint.
The variable $\ph$ stands for an order parameter
and  is
usually taken between $-1$ and $1$\juerg{; it represents} the healthy cell case 
and the tumor phase, respectively. Moreover, $\m$ indicates the chemical potential 
for $\ph$, whereas $\s$ denotes the nutrient extra-cellular
water concentration. This latter quantity is usually normalized between 
$0$ and $1$, conveying that these values model the nutrient-poor and
the nutrient-rich cases. The symbols $\a$ and $\b$ represent positive
constants; \lastrev{let us just note that the term $\b\dt\ph$ in the second equation
corresponds to the classical term
of the viscous Cahn--Hilliard equation, while 
the term $\a\dt\m$ gives to equation \eqref{State_first} a parabolic structure
with respect to $\m$. For more details on these relaxation terms, let us 
refer to \cite{CGH,CGRS_VAN,CGRS_ASY}.}
The capital letters $A,B,D,P,\ch$ denote
positive coefficients that stand for
the apoptosis rate, nutrient supply rate, nutrient consumption rate, proliferation
rate, and chemotaxis coefficient, in \juerg{this} order.
\lastrev{In addition, let us point out that the contributions $\ch \s$ and $\ch\Delta \ph$
model \juerg{pure chemotaxis, namely,} the movement of tumor
cells towards regions of high nutrients, and the active transport 
that describes the movement of the nutrient towards the
tumor cells (see \cite{GARL_1,GARL_2,GAR} for more details).
}
Furthermore, the function $h$ \last{has been considered
in \cite{GLSS} (see also \cite{GARL_2, GARL_1, GARL_4, CLLW, WLFC})} 
as an interpolation function between
$-1$ and $1$ in order to have $h(-1)=0$ and $h(1)=1$, so that
the mechanisms modelled by the terms
$(P\s -A - u) h(\ph)$ and $D \s h(\ph)$
are switched off in the healthy case, which
corresponds to $\ph=-1$, and are fully active in the tumorous case $\ph=1$.
Besides, the term $\ss$ stands for a nonnegative constant modelling the nutrient 
concentration in a pre-existing vasculature. 
For further details on the model,
we refer the reader to \cite{GLSS} (see also \cite{GARLR,SM}).
Lastly, the term $F'$ is the derivative of a 
double-well nonlinearity. Typical examples for this 
nonlinearity are the regular potential
\begin{align}
  F_{reg}(r) = \frac14(r^2-1)^2 \, 
  \quad \hbox{for } r \in \erre,
  \label{F_reg}
\end{align}
and\juerg{, more relevant for applications, the logarithmic potential}
\begin{align}
   F_{log}(r) = (1+r)\ln (1+r)+(1-r)\ln (1-r) - k r^2 \,
  \quad \hbox{for } r \in (-1,1),
  \label{F_log}
\end{align}
where \juerg{$k>1$ so that $F_{log}$ is nonconvex}.
Eventually, the terms $u$ and
$w$ are \juerg{source} terms acting as control variables.
It is worth noting that we are considering two control variables: $u$
in the phase equation and $w$ in the nutrient equation.
In the previous contributions \cite{CGRS_OPT,S,S_a,S_b,CRW}, the control variable was
placed in the nutrient equation,
so that it designs an external medication or some nutrient supply.
On the other hand, different authors consider the control variable in the phase equation
\juerg{(see, e.g., \cite{GARLR,EK,EK_ADV})}, multiplied by $h(\ph)$ in order to have the 
action of the control 
only in the meaningful region. In that case,
it models the introduction of cytotoxic drugs into the system, which 
has the purpose of eliminating
the tumor cells.
Thus, with our choice we include both \juerg{these} cases in this paper.

We are now in a position to introduce the distributed optimal control 
problem we are going to deal with. It consists of \juerg{finding} a solution 
to the following minimization problem:

\vspace{2mm}\noindent
$\CP$ Minimize the tracking-type cost functional
\begin{align}
	\non
	\J(\ph,\s,u,w) : = &
		\frac \bO 2 \iO |\ph(T)-\phO|^2		
		+ \frac \bQ 2 \intQ |\ph-\phQ|^2
		+ \frac \bG 2 \iO |\s(T)-\sO|^2		
		\\ & \qquad 
		+ \frac \bS 2 \intQ |\s-\sQ|^2
		+ \frac \bu 2 \intQ |u|^2
		+ \frac \bw 2 \intQ |w|^2
		\label{cost}
\end{align}
\juerg{subject to the condition} that $(\m,\ph,\s)$ solves the {\em state system}
\Stsys\ \juerg{for a control pair $(u,w)$ belonging} to the \lastjuerg{control box}
\begin{align}
	\label{Uad}
	\Uad:= \graffe{(u,w) \in (L^\infty (Q))^2: 
		u_* \leq u \leq u^*\ \aeQ, \ w_* \leq w \leq w^*\ \aeQ},
\end{align}
where $u_*,u^*,w_*$ and $w^*$ denote some prescribed functions in $L^\infty(Q)$.
Moreover, let us point out that the physical meaning of the term $u h(\ph)$ 
in the state system requires the control $u$ to be nonnegative. 
Hence, in the following we will \juerg{always assume that the lower bound satisfies
$u_*\ge 0$ almost everywhere in $Q$}.

As far as control problems for tumor growth models
are concerned, the contributions are \juerg{still scarce. To} our knowledge,
the first optimal control problem
governed by a tumor growth model similar to the one given above is \cite{CGRS_OPT}.
There, the control problem was investigated for the case of
regular potentials enjoying polynomial growth.
Then, by adding two relaxation terms, a similar optimal control
problem was tackled in \cite{S} by extending the generality
of the potentials by allowing singular, \juerg{but still smooth}, 
potentials like the logarithmic potential to be considered.
Next, the same author, using the \socal\ deep quench asymptotic technique,
proved in \cite{S_DQ} how nonsmooth potentials like the 
double obstacle \juerg{potential can also be admitted}.
Then, exploiting the results known for the case $\a,\b>0$,
in the following works \cite{S_a,S_b} the author \juerg{showed} that it is possible to 
let $\a$ and $\b$ \juerg{approach} zero separately \juerg{in order to recover}
 the existence of optimal controls
and \juerg{to characterize} the corresponding 
first-order necessary conditions for optimality.
We also refer to \cite{GARLR}, where
an optimal treatment time has been performed for a similar system, namely for
system \Stsys\ with the choices $\ch=\a=\b$ and $w\equiv 0$; see also \cite{CRW}, where
a similar control problem was investigated for a different model.
Moreover, let us mention \cite{SW}, where an optimal control problem
for the two-dimensional Cahn--Hilliard--Darcy system
with mass sources is addressed.
We also point out \cite{EK_ADV,EK}, where the optimal control for a
Cahn--Hilliard--Brinkman type system has been tackled.
Lastly, we refer to \cite{SM}, where a different kind of 
control problem, known as sliding
mode control, was performed for a system that is very close to \Stsys.

We now comment on \Stsys. Let us emphasize that,
once the well-posedness of the state system is established,
we can properly define the \juerg{{\em control-to-state operator} that
assigns to a given control $(u,w)$ the 
unique corresponding solution to \Stsys, }
\begin{align}
	\label{controltostate}
	\S:(u,w)\mapsto \S(u,w):=\bigl(\m,\ph,\s),
\end{align}
and attains values in a proper Banach space.
Then, we are in a position to eliminate the state variable appearing in
the cost functional \eqref{cost} by
expressing them as functions of the control. This leads to \juerg{the {\em reduced cost functional}}
\begin{align}
\label{Jred}
	\Jred(u,w):= \J (\S_2(u,w),\S_3(u,w),u,w),
\end{align}
where $\S_2(u,w)$ and $\S_3(u,w)$ denote the second and third component of $\S$, respectively.
At this formal stage, let us point \lastjuerg{out} that from standard 
results of convex analysis (see, e.g., \cite{Trol,Lions_OPT}) \juerg{it follows}
the formal first-order necessary condition for 
optimality characterized by the variational inequality
\begin{align}
	\label{formal_foc}
	D\Jred(\uopt,\wopt)(u-\uopt,w-\wopt) 
	\geq 0 \quad \hbox{for every $(u,w) \in \Uad$,}
\end{align}
where $D\Jred $ stands for the derivative of the reduced cost functional
in a proper \last{functional analytic} sense.

Therefore, summing up, in this contribution we aim
at solving the constrained minimization problem
\begin{align}
	\non
	\CP \quad 
	& \hbox{Minimize ${\J}(\ph,\s,u,w)$}
	\hbox{ subject to the control contraints \eqref{Uad} and under the}
	\\ & \non
	\hbox{requirement that the variables $(\ph, \s)$ 
	yield a solution to \Stsys,}
\end{align}
and pointing out the corresponding first-order necessary continuations
for optimality.

The paper is organized as follows.
The next section brings the mathematical 
framework and gathers the obtained results.
From Section \ref{SEC_STATE} onward, we proceed with 
the proof of the statements.
The well-posedness and the continuous dependence results for
the state system \Stsys\ are addressed in Section~\ref{SEC_STATE}, 
while Section~\ref{SEC_CONTROL_PROBLEM} is completely devoted
to the corresponding control problem. Namely, we prove in this last section
the existence of optimal controls and 
\juerg{derive} the corresponding first-order necessary conditions for optimality.

\section{Mathematical setting and main results}
\label{SEC_MATHEMATICAL_SETTING}
\setcounter{equation}{0}
To begin with, let us point out some notation.
\juerg{As far as} the functional spaces are concerned, it is 
convenient to set
\begin{align}
	\non
	H:= \Ldue, \quad V:= \Huno, \quad W:=\graffe{v \in \Hdue : \dn v = 0 \hbox{ on } \Gamma},
\end{align}
and \juerg{to} endow $H,V,W$ with their standard norms.
Furthermore, for an arbitrary Banach space $X$,
we denote by $\|\,\cdot\,\|_X$ its norm, $X^*$ its topological dual, 
\juerg{and by} $\langle\, \cdot\,,\, \cdot\rangle\juerg{_X} $ the duality product between $X^*$ and $X$.
Likewise, for every $1\leq p \leq \infty$, we use the symbol $\|\,\cdot\,\|_{p}$ to indicate the usual norm 
in $L^p(\Omega)$. \juerg{Notice} that $(V,H,\Vp)$ forms a Hilbert triple,
that is, the injections $V \subset H \equiv H^* \subset V^*$ 
are both continuous and dense, \juerg{where we have the identification}
\begin{align*}
	\< u,v >\juerg{_V} = \iO uv
	\quad \hbox{for every $u \in H$ and $v \in V$}.
\end{align*}
Furthermore, it is convenient to denote the parabolic cylinder and its boundary by
\begin{align}
	\non
	Q_{t}&:=\Omega \times (0,t) \quad \hbox{and} \quad \Sigma_{t}:=\Gamma \times (0,t)
	\quad \hbox{for every $t \in (0,T]$,}
	\\ \non
	Q&:=Q_{T}, \quad \hbox{and} \quad \Sigma:=\Sigma_{T}.
\end{align}
%

For the potential $F$, we generally assume:

\vspace{2mm}\noindent
{\bf (F1)} \,\,$F=F_1+F_2$, where $F_1:\erre \to [0,+\infty]$ is convex and lower
semicontinuous \linebreak
\hspace*{12.2mm}with $F_1(0)=0$\last{.} 

\vspace{2mm}\noindent {\bf (F2)} \,\,\last{There exists an interval $(\rmin,\rmax)$, with $-\infty \leq r_- <0 < r_+ \leq + \infty$, such that}
\linebreak
\hspace*{12.5mm}\last{the restriction 
of $F_1$ to $(\rmin,\rmax)$ is differentiable 
with derivative
$F'_1$.}

	
\vspace{2mm}\noindent
{\bf (F3)} \,\,
	\last{$F_2\in C^3(\erre)$}, and $F_2'$ is \Lip\ continuous with \Lip\ constant $L>0$.

\vspace{2mm}\noindent
{\bf (F4)} \,\,	$F_{| _{\last{(\rmin,\rmax)}}} \in C^3(\rmin,\rmax),
	\hbox{ and } \lim\limits_{r \to r_{\pm}} F'(r) = \pm \infty$.

\noindent It is worth noting that both \eqref{F_reg} and \eqref{F_log} do fit
the above framework with the choices
$(\rmin,\rmax)=(-\infty,+\infty)$ and $(\rmin,\rmax)=(-1,1)$, respectively, 
so that
they are allowed to be considered.

For the initial data introduced above, we \juerg{make the following assumptions:}

\vspace{2mm}\noindent
{\bf (A1)} \,\,	$\ph_0 \in W,\,\,\, \m_0 \in \Huno\cap \Lx\infty ,\,\,\, \s_0 \in \Huno\cap \Lx\infty.$
	
\vspace{2mm}\noindent
{\bf (A2)} \,\,	$\rmin < \inf \ph_0 \leq \sup \ph_0 < \rmax,$ \,\,\,
whence $F(\ph_0),\,\,F'(\ph_0)\in\Lx\infty.$

\last{
\noindent
Notice that the above requirement
can be restrictive for the case of singular potentials. 
For instance, in the case of the logarithmic potential, we have
$r_\pm=\pm 1$ so that \last{by {\bf (A2)} the pure phases (tumor and
healthy tissue) can only be approximated by the initial datum.}}

\vspace{2mm}\noindent
For the other appearing constants and target functions, we
postulate:

\vspace{2mm}\noindent
{\bf (A3)} \,\,$h\in \last{C^2(\erre)\cap W^{2,\infty} (\erre)}$,
 and \last{$ h $ is positive on $(\rmin,\rmax)$.}

\vspace{2mm}\noindent
{\bf (A4)} \,\, $\a, \b, \ch $ are positive constants, while $P, A,B,D, \s_s$ are nonnegative constants.	

\vspace{2mm}\noindent
{\bf (A5)} \,\,	$\bO,\bQ, \bG, \bS, \bu,\bw $  \,are nonnegative constants, but not all zero.
	
\vspace{2mm}\noindent
{\bf (A6)} \,\, 	$\phQ, \sQ \in L^2(Q),\,\,\,\phO, \sO \in \Lx2.$

\vspace{2mm}\noindent
\last{Note that {\bf (A3)} entails that $h$, $h'$ are Lipschitz continuous in $\erre$. 
Let us denote by $\hinf,\hinf'$ the upper 
bounds for the $L^\infty (\erre)$ norms of $h$ and $h'$, 
and by $L_h$ the \Lip\ constant of $h$, respectively.
Let us also point out that as $L_h$
we can simply take ${h'}\!_\infty$.}

\noindent Moreover, we assume that the control box $\Uad$ is defined by \eqref{Uad}, and that

\vspace{2mm}\noindent
{\bf (A7)} \,\,
	$u_*,u^*\in L^\infty(Q) \hbox{ with }\,  0\leq u_* \leq u^* \,\, 
	\,\aeQ, \ w_*,w^*\in L^\infty(Q) \hbox{ with }
		 w_*\leq w^* \,\,\,\mbox{a.e.}$ \linebreak \hspace*{13.9mm}in \,$Q$.

\vspace{2mm}\noindent
The latter condition implies that $\Uad$ is a closed and convex \lastjuerg{subset} of $L^2(Q)$.
On the other hand, it will be sometimes convenient to work with  
an open \juerg{superset of $\Uad$. We therefore fix some constant $R>0$ such that \last{the open ball}}
\begin{align}
\label{defUR}
	&{ \mbox{${\UR := }\{(u,w)\in L^\infty(Q)\times L^\infty(Q):\,
	\norma {(u,w)}_{L^\infty(Q)\times L^\infty(Q)} < R \} $\, contains $\,\Uad$.}}
\end{align}

\vspace{1mm}
\Brem
Before diving into the well-posedness result, let us point
out a classical issue of control theory. The well-posedness result 
\juerg{to be presented below} is given in a rather strong setting; 
this is motivated by the control problem under investigation. \juerg{However,
system \Stsys\ can be provided with} a notion 
of weak solutions in a rather mild setting.
Moreover, it is also possible to extend the analysis
for the potentials \last{taking} into account singular and nonregular potentials like
the well-known double obstacle potential. For this, a Yosida regularization
of the maximal monotone operator $F_1'$ has to be introduced.
Clearly, the pointwise formulation \Stsys\ has then to be replaced by 
a suitable variational formulation. \juerg{Let us just sketch the expected result here:
provided we assume $\m_0,\ph_0,\s_0\in \Lx2$ for the initial data and a potential
that fulfills {\bf (F1)}--{\bf (F3)}, we can prove existence}
and uniqueness of a weak solution  such that
$\m,\ph,\s \in \H1 \Vp \cap \L\infty H \cap \L2 V$. Note that uniqueness 
will follow from the first continuous dependence estimate that we perform below 
(cf. \eqref{cont_dep_first}), which perfectly complies with the above notion of 
weak solutions.		
\Erem
First, let us present the result regarding the existence and
uniqueness of a strong solution to the system \Stsys.
\begin{theorem}
	\label{THM_well-posedness_state}
	Assume that \,{\bf (F1)}--{\bf (F4)}, {\bf (A1)}--{\bf (A4)}, \juerg{and {\bf (A7)}, are fulfilled and 
	that} $(u,w)\in \UR$.
	Then the state system \Stsys~admits a unique
	solution $(\m,\ph,\s)$ \juerg{with the} regularity
	\begin{align}
		\label{reg_state_first}
		\ph & \in W^{1,\infty}(0,T;H) \cap \H1 V \cap \L\infty W\,,
		\\ 
		\label{reg_state_second}
		\m,\s & \in \H1 H \cap \L\infty V \cap \L2 W \cap L^\infty (Q)\,.
	\end{align}
	\Accorpa\regstate reg_state_first reg_state_second
	Moreover, there exists a positive constant $K_1$, \juerg{which depends only on 
	$\Omega$, $T$, $R$, $\a$, $\b$, and the data of the system}, such that
	\begin{align}
	&\non
	\norma{\ph}_{W^{1,\infty}(0,T;H) \cap \H1 V \cap \L\infty W}
	+ 	\norma{\m}_{\H1 H \cap \L\infty V \cap \L2 W \cap L^\infty (Q)}
	\\ & \quad 
	+ \norma{\s}_{\H1 H \cap \L\infty V \cap \L2 W \cap L^\infty (Q)}
	\leq K_1\,.
	\label{stima_reg}
	\end{align}
	In addition, there exist some constants $r_{*}$ and $r^{*}$, which
	satisfy $ \rmin < r_{*}\leq r^{*} < \rmax $ and depend only on the data
	of the system, such that
\begin{align}
	\label{separation_result}
	r_* \leq {\ph} \leq  r^* \quad \aeQ\,.
\end{align}
Finally, there exists a positive constant $K_2$, \juerg{which depends only on $\Omega$, 
$T$, $R$, $\a$, $\b$, and the data of the system}, such that
\begin{align}
	\label{stima_postsep}
	\norma{\ph}_{L^\infty(Q)} + \max_{i=1,2,3} \norma{F^{(i)}(\ph)}_{L^\infty(Q)}
	\leq
	K_2\,.
\end{align}
\end{theorem}
\last{
\noindent
Let us point out that \eqref{separation_result} turns out to be significant
for the case of singular potentials such as the logarithmic potential.
In that situation, it guarantees that, as soon as \last{\bf (A2)} is satisfied, the phase variable
stays away from the pure \last{phases uniformly} during the evolution.
This fact, known as the separation principle, in turn, entails that
\eqref{stima_postsep} holds.}
\vspace{1mm}
\Bthm
\label{THM_cont_dep}
Suppose that {\bf (F1)}--{\bf (F4)} and {\bf (A1)}--{\bf (A4)} are fulfilled.
\juerg{Then there exists a positive constant $K_3$, which depends only on $\Omega$, 
$T$, $R$, $\a$, $\b$, and the data of the system, such that the following holds true:
whenever two control pairs $(u_i,w_i) \in \UR$, $i=1,2$, are given  and $(\m_i,\ph_i,\s_i)$,
$i=1,2$, are the \andrealast{corresponding} states, then }
\begin{align}
	\non &
	\norma{\m_1-\m_2}_{\H1 H \cap \L\infty V \cap \L2 W}
	+ \norma{\ph_1-\ph_2}_{\H1 H \cap \L\infty V \cap \L2 W}
	\\ & \qquad  \non
	+ \norma{\s_1-\s_2}_{\H1 H \cap \L\infty V \cap \L2 W}
	\\ & \quad
	\leq 
	K_3 \, \bigl(\norma{u_1-u_2}_{\L2 H}+\norma{w_1-w_2}_{\L2 H}\bigr)\,.
	\label{cont_dep_final}
\end{align}
\Ethm

\vspace{1mm}
\juerg{For the optimal control problem $\CP$, we will show the following existence result:}
\Bthm
\label{THM_existence_opt_control}
Assume that {\bf (F1)}--{\bf (F4)} and {\bf (A1)}--{\bf (A7)}  are satisfied.
Then the control problem $\CP$ admits at least \juerg{one} solution.
\Ethm

\vspace{2mm}\juerg{Finally, we formulate the first-order necessary optimality conditions for
$\CP$ that will \pier{be} shown below. To this end, we} assume that $(\uopt,\wopt)$ and $(\bm,\bph,\bs)$ stand for
some fixed control and its associated state, respectively.
Sometimes, the same notation is employed to refer to an optimal
control with the corresponding optimal state; anyhow,  we will specify \juerg{this}
whenever it is the case. \juerg{In the course of our analysis, it will be necessary to establish the
Fr\'echet differentiability of the control-to-state operator ${\cal S}:(u,w)\mapsto(\mu,\phi,\sigma)$
in suitable Banach spaces. To this end, the unique solvability 
of the corresponding linearized system will have to be shown. This system has 
for every pair $(k,l)\in (L^2(Q))^2$ the following form: }
\begin{align}
	\label{Lin_first}
	&\a \dt \et
	+ \dt \xi- \Delta \et
	=
	(P\z - k) h(\bph)
	+( P\bs - A -\uopt) h'(\bph)\xi
	\quad &&\hbox{in $\, Q$},
	\\
	\label{Lin_second}
	&\et
	=
	\b\dt \xi
	-\Delta \xi
	+ F''(\bph)\xi
	- \ch \z
 	\quad &&\hbox{in $\,Q$},
	\\
	\label{Lin_third}
	&\dt \z
	- \Delta \z
	+B\z
	=
	\andrea{-\ch \Delta \xi}
	- D \z h(\bph)
	- D \bs h'(\bph)\xi
	+ l
	\quad &&\hbox{in $\,Q$},
	\\
	&\dn \et=\dn \xi=\dn \z =0
 	\quad &&\hbox{on $\,\Sigma$},
	\label{Lin_BC}
	\\ 
	&\et(0)=\xi(0)=\z(0)=0
  	\quad &&\hbox{in $\,\Omega.$}
	\label{Lin_IC}
\end{align}
\Accorpa\Linsys Lin_first Lin_IC
Here, the well-posedness result follows.
\begin{theorem}
	\label{THM_well-posedness_lin}
	Assume that {\bf (F1)}--{\bf (F4)}, {\bf (A1)}--{\bf (A4)}, and {\bf (A7)}, are satisfied. Then
		the linearized system \Linsys~admits for every $(k,l)\in (L^2(Q))^2$ a unique
	solution $(\et,\xi,\z)$ with the regularity
	\begin{align}
	\label{reg_lin}
		\et,\xi,\z & \in \H1 H \cap \L\infty V \cap \L2 W.
	\end{align}
\end{theorem}
\juerg{Notice that Theorem~\ref{THM_cont_dep}
also entails the \Lip\ continuity of the 
control-to-state operator $\S$ between suitable Banach spaces.
We even have Fr\'echet differentiability, as the following result states.} 

\Bthm
\label{THM_frechet}
Assume that {\bf (F1)}--{\bf (F4)}, and {\bf (A1)}--{\bf (A4)} are satisfied,
\juerg{and let $(\uopt,\wopt)\in \UR$ be} a fixed control 
with the corresponding state $(\bm,\bph,\bs)$.
Then the control-to-state operator $\S$ is \Frechet\ differentiable at
$\juerg{(\uopt,\wopt)}$ as a mapping from $(L^\infty(Q))^2$ into the Banach space $\cal Y$,
where
\begin{align}
	\non
	{\cal Y} :=& \bigl( \C0 H \cap \L2 V\bigr) \times 
	\bigl( \H1 H \cap \L\infty V  \bigr) 
	\\ & 
	\times
	\bigl( \C0 H \cap \L2 V \bigr).
	\label{pier2}
\end{align}
Moreover, for every $(k,l)\in (L^\infty(Q))^2$ the derivative of $\,\S\,$ at $\,\juerg{(\uopt,\wopt)}\,$
is given by the identity $[D\S (\uopt,\wopt)](k,l)=(\et,\xi,\z)$, where $(\et,\xi,\z)$
is the unique solution to the linearized system \Linsys\ corresponding to $(k,l)$.
\Ethm

\Bthm
\label{THM_foc_prima}
Assume that {\bf (F1)}--{\bf (F4)} and {\bf (A1)}--{\bf (A7)} are fulfilled, \juerg{and let 
$(\uopt,\wopt)$ be an optimal control with associated state $(\bm,\bph,\bs)$.
Then it holds that}
\begin{align}
	& \non
	\bO \iO (\bph(T)-\phO)\xi(T)
	+ \bQ \intQ (\bph-\phQ)\xi
	+ \bG \iO (\s(T)-\sO)\z(T)
	+ \bS \intQ (\bs-\sQ)\z
	\\ & \quad
	+ \bu \intQ \uopt(u-\uopt)
	+ \bw \intQ \wopt(w-\wopt)
	\geq 0 \quad \hbox{for every $(u,w)\in \Uad$},
	\label{foc_first}
\end{align}
where the triple $(\et,\xi,\z)$ is the unique solution to the linearized system \Linsys\ corresponding to $k=u-\uopt$ and $ l= w-\wopt$,
respectively.
\Ethm

Analyzing the above variational inequality, one realizes that it 
is not very useful in numerical applications, since for every
possible step of the approximation
one has to solve the state system and also the linearized system in order
to have $\xi$ and $\z$ at disposal.
For this reason, a classical tool is to introduce the \socal~adjoint system
in order to eliminate these variables. In fact, provided that we choose
this auxiliary system properly, the linearized variables can be eliminated from
\eqref{foc_first}. The adjoint system to \Stsys\ can be obtained 
by the formal Lagrangian method described, e.g., in \cite{Trol}, using
integration by parts. \juerg{Following this route, we arrive at  the following (formal) version
of the adjoint system:} 
\begin{align}
	\label{Adj_first}
	& 	-\a \dt q 	- \Delta q 	- p 	\,=\,	0
	\qquad \hbox{in $\, Q$},
	\\
	& \non
	-\dt q 
	-\b \dt p 
	- \Delta p 
	+\ch \Delta r
	+F''(\bph)p
	- (P\bs -A -\uopt)h'(\bph) q 
	+ D\bs h'(\bph) r
	\\
&	\hspace*{6.7mm}=\, 	\bQ (\bph-\phQ)
 	\qquad \hbox{in $\,Q$},
 	\label{Adj_second}
	\\
	\label{Adj_third}
	& -\dt r
	- \Delta r
	+ Br
	+ D h(\bph)r
	- \ch p
	- P h(\bph)q
	\,=\, \bS (\bs - \sQ)
	\qquad \hbox{in $\,Q$},
	\\
	&\dn q\,=\,\dn p\,=\,\dn r \,=\,0
 	\qquad \hbox{on $\,\Sigma$},
	\label{Adj_BC}
	\\ 
	& q(T)\,=\,0, \ \b p(T)\,=\, \bO (\bph(T)-\phO), \ r(T)\,=\,\bG(\bs(T)-\sO)
  	\qquad \hbox{in $\,\Omega.$}
	\label{Adj_IC}
\end{align}
\Accorpa\Adjsys Adj_first Adj_IC

\juerg{Observe that this is a backward-in-time system with
final conditions belonging to $\Lx2$ (see assumption {\bf (A6)}), so that}
we cannot expect to recover a strong solution\andrealast{.}
Therefore, instead of considering 
the pointwise equations \accorpa{Adj_second}{Adj_third}, we note that
the variables $p$ and $r$ should be understood as  weak solutions in the
following sense:
\begin{align}
\label{pweak}
&-\iO\dt q\,v-\langle\beta\dt p,v\rangle\juerg{_V}+\iO\nabla p\cdot\nabla v
-\chi\iO\nabla r\cdot\nabla v + \iO F''(\overline{\varphi})\,p\,v
+\iO D\,\overline{\sigma}\,h'(\overline{\varphi})\,r\,v \nonumber\\
&-\iO(P\overline{\sigma}-A-\overline{u})h'(\overline{\varphi})q\,v
\,=\iO\gamma_2(\overline\varphi-\phi_Q)\,v\quad\,\mbox{for all $v\in V$ and a.e. in $(0,T)$},\\
\label{rweak}
&	-\<\dt r,v>\juerg{_V}
	+\iO \nabla r\cdot \nabla v
	+\iO Br v
	+\iO D h(\bph)r\, v
	-\iO \ch p \,v
	-\iO P h(\bph)q\, v
	\nonumber\\ 
	&  
	= \iO\bS (\bs - \sQ) v
	\quad \,\hbox{for all $v\in V$  and a.e. in $(0,T)$},
\end{align}
where, for simplicity, we avoided writing the time variable explicitly.
We have  the following well-posedness result.
\Bthm
\label{THM_wellposed_adjoint}
Assume that {\bf (F1)}--{\bf (F4)} and {\bf (A1)}--{\bf (A7)} are fulfilled, \juerg{and let 
$(\uopt,\wopt)$ be an optimal control with associated state $(\bm,\bph,\bs)$.
Then the adjoint system \Adjsys\ has} 
a unique solution such that
\begin{align}
	\label{reg_adj_prima}
	p,r &\in \H1 \Vp \cap \C0 H \cap \L2 V,
	\\
	q &\in \H1 H \cap \C0 V \cap \L2 W.
	\label{reg_adj_seconda}
\end{align}

\Ethm

\Bthm
\label{THM_foc_final}
Assume that {\bf (F1)}--{\bf (F4)} and {\bf (A1)}--{\bf (A7)} are fulfilled, \juerg{and let 
$(\uopt,\wopt)$ be an optimal control with associated state $(\bm,\bph,\bs)$ and adjoint state
$(p,q,r)$. Then it holds the variational inequality}
\begin{align}
	& 
	\label{foc_final}
	\intQ (-h(\bph) q + \bu \uopt)(u-\uopt)
	+  \intQ (r + \bw\wopt)(w-\wopt)
	\geq 0 \quad \hbox{for every $(u,w)\in \Uad$}.
\end{align}
Moreover, whenever $\bu\not = 0$, then $\uopt$ is nothing but the 
$\L2 H$-orthogonal projection of $\bu^{-1} h(\bph) q\,$ onto the closed and convex set 
\juerg{$\{u\in L^2(Q):\,u_*\le u\le u^* \,\mbox{ a.e. in } \,Q\}$}.
Likewise, if $\bw\not = 0$, then $\wopt$ reduces to 
the $\L2 H$-orthogonal projection of $-\bw^{-1} r$ onto 
\juerg{$\{\pier{w}\in L^2(Q):\,w_*\le w\le w^* \,\mbox{ a.e. in } \,Q\}$}.
\Ethm

Furthermore, since $\Uad$ is actually a control box, it is possible
to explicitly characterize the projection and obtain a pointwise condition.
\Bcor
Let {\bf (F1)}--{\bf (F4)} and {\bf (A1)}--{\bf (A7)} be fulfilled, and let $\bu>0$. 
Then, the optimal control \juerg{component} \,$\uopt\,$ is implicitly characterized by
\begin{align}
	\non
	\uopt(x,t)=\max \bigl\{ 
	u_*(x,t), \min\graffe{u^*(x,t),\bu^{-1} h(\bph(x,t)) q(x,t)} 
	\bigr\} \quad \aaQ .
\end{align}
Likewise, if $\bw >0$, then
\begin{align}
	\non
	\wopt(x,t)=\max \bigl\{ 
	w_*(x,t), \min\graffe{w^*(x,t),-\bw^{-1} r(x,t)} 
	\bigr\} \quad \aaQ .
\end{align}
\Ecor

Let us emphasize a consequence which is of \sfw~importance for the numerical approach.
Comparing the expected theoretical condition \eqref{formal_foc}
with the explicit condition \eqref{foc_final}, via Riesz's representation
theorem, the gradient of the reduced cost functional
can be recovered as $\nabla {\Jred}(\uopt,\wopt)=(\last{-}h(\bph) q + \bu \uopt\,,\,r + \bw\wopt)$.
Hence, for the numerical approach, the optimal control problem can be viewed
as a constrained minimization of a function, $\J_{red}$, whose gradient is known
(think of the \wk~projected conjugate gradient method).

In the remainder of this section, we recollect some \wk\ \lastjuerg{results
that} will prove useful later on.
To begin with, we recall 
the standard Sobolev continuous embedding
\begin{align}
	\label{VinL6}
	\Huno \hookrightarrow L^q(\Omega)
	\quad \hbox{for every $q\in [1,6]$.}
\end{align}
Furthermore, we often make use of Young's inequality
\begin{align}
  ab \leq \delta a^2 + \frac 1 {4\delta} \, b^2
  \quad \hbox{for every $a,b\geq 0$ and $\delta>0$}.
  \label{young}
\end{align}

As far as the constants are concerned, let us set 
our convention once and for all: the symbol small-case $c$ is used to indicate every constant
that depends only on the structural data of the problem, such as
$T$, $\Omega$, $R$, $\a$ or $\b$, the shape of the 
nonlinearities, and the norms of the involved functions. 
On the other hand, with capital letters 
we specify particular constants
that will be referred to later on.
Therefore, the meaning of the constant $c$ may change from line to line.

\section{The state system}
\label{SEC_STATE}
\setcounter{equation}{0}
\subsection{Well-posedness of the state system}

\Bdim[Proof of Theorem~\ref{THM_well-posedness_state}]
Here, we \last{do not provide all the details}, since the approach is quite standard.
Anyhow, let us point out that the argument can be made rigorous by \last{using} an approximation technique,
e.g., within a Faedo--Galerkin scheme along with the introduction of the Yosida
approximation for $F_1'$. In fact, since the framework for the potential
settings is rather general, we cannot assume $F_1'$ to
be \Lip\ continuous, in general. 
\last{
In this direction, let $\eps \in (0,1)$ and, for every $r \in \erre$, let us set
\begin{align*}
	F_{1,\eps} (r):= 
		\min_{s\in\erre} \Big( \frac 1 {2\eps} (s-r)^2 + F_1(s)\Big),
		\quad
		F'_{1,\eps}(r):= \frac d {dr} F_{1,\eps} (r),
		\quad
		F_{\eps}(r):= F_{1,\eps}(r)+F_2(r).
\end{align*}
It turns out that $F_{1,\eps} \in C^1(\erre)$ 
and that the Yosida regularization $F'_{1,\eps}$ is \Lip\ \last{continuous} 
(see, e.g., \cite[Prop.~2.11, p.~39]{BRZ}).
Furthermore, \last{for every $r \in \erre$ the following properties
\begin{align}
	0 \leq F_{1,\eps}(r) \leq F_{1}(r),
	\quad
	F_{1,\eps}(r) \nearrow F_{1}(r) \quad
	\hbox{monotonically as $\eps \searrow 0$}
\label{pier3}
\end{align}
hold, as well as (cf.~\cite[Prop.~2.6, p.~28]{BRZ})
\begin{align}
	\hbox{for all $r\in  (r_-,r_+)$} \quad
	|F'_{1,\eps}(r) |\nearrow | F_{1}'(r)|  \quad
	\hbox{monotonically as $\eps \searrow 0$}.
\label{pier4}
\end{align}}%
Hence, the idea is as follows:
first we aim at discussing the well-posedness of the \last{approximation of} \Stsys, 
namely \Stsys\ \last{with $F'$ replaced by} $F'_{\eps}$, and
then, on account of a priori estimates and monotonicity arguments,
using this result to ensure the existence of a solution for the original system.
For the sake of simplicity, we denote by $(\m_\eps, \ph_\eps, \s_\eps)$
the solution to the approximated system but
we avoid writing the subscript $\eps$ in the calculations below.
Only at the end of each calculation the correct notation is employed.}
The proof of the uniqueness will follow 
as a direct consequence of Theorem~\ref{THM_cont_dep}.

\noindent {\bf First estimate:} \,\,\,
To begin with, we add to both sides of \eqref{State_second} the term $\ph$.
Then, we multiply \eqref{State_first} by
$\m$, the new \eqref{State_second} by $\dt\ph$, \juerg{ \eqref{State_third} by $\s$, and add
the resulting identities}. 
Next, we integrate over $Q_t$, for an arbitrary $t\in(0,T]$, and by parts.
After a cancellation of terms and some rearrangements, we infer that
\begin{align*}
	& \frac \a2 \normaH{\m(t)}^2
	+ \I2 {\nabla \m}
	+ \b \I2 {\dt \ph}
	+ \frac 12 \normaH{\ph(t)}^2
	+ \frac 12 \normaH{\nabla \ph(t)}^2
	\\ & \qquad
	\last{{}+ \iO F_{1,\eps}(\ph(t)) }
	+ \frac 12 \normaH{\s(t)}^2
	+ B \intQt |\s|^2 
	+  \I2 {\nabla \s}
	\\ &\quad
	= 
	\frac \a2 \normaH{\m_0}^2
	+ \frac 12 \normaH{\ph_0}^2
	+ \frac 12 \normaH{\nabla \ph_0}^2
	+ \frac 12 \normaH{\s_0}^2
	\last{{}+  \iO F_{1,\eps}(\ph_0)}
	\\ & \qquad
	+ \intQt (P\s -A- u) h(\ph) \m
	+ \ch  \intQt \s \dt \ph
	+\intQt (\ph -F_2'(\ph))\dt \ph
	\\& \qquad
	+ \ch \intQt \nabla \ph \cdot\nabla \s
	+ \intQt B \ss \s
	- \intQt D h(\ph) |\s|^2
	+ \intQt w \s.
\end{align*}
\juerg{Obviously, all of the summands on the \lhs\ are nonnegative, 
and the first \last{four} summands on the \rhs\ are
bounded, by virtue of {\bf (A1)}, {\bf (A2)}, and the general assumptions on $F_1$ and $F_2$. 
\last{
Besides, \last{\eqref{pier3}} implies that the fifth term verifies 
\begin{align*}
	\iO F_{1,\eps}(\ph_0) \leq
	\iO F_{1}(\ph_0) \leq c.
\end{align*}
}%
It remains to estimate 
the remaining terms on the \rhs, which we denote by $I_1,...,I_7$, in this order. This can easily be done by
means of Young's inequality.}
In fact, we have that
\begin{align*}
	|\juerg{I_1}|
	& \leq
	\I2 \s
	+ T|\Omega|
	+ \I2 u
	+ \frac {\hinf^2 \bigl( P^2 + A^2 +1 \bigr)}4 \I2 {\m}.
\end{align*}
Furthermore, we also infer that
\lastrev{
\begin{align*}
	\juerg{\sum_{i=2}^{7}}|I_i|
	&\leq
	\frac \b2 \I2 {\dt\ph}
	+ \frac {\ch^2}\b \I2 \s
	+ \frac {2(1+L^2)}\b \I2 \ph
	+ \frac {\ch^2}2 \I2 {\nabla \ph}
	\\ & \quad
	+ \frac 12 \I2 {\nabla \s}
	+ \frac 12 \I2 \s
	+ \frac {B^2 \ss^2}2 \, T\,|\Omega|
	+ D \hinf \I2 \s
	+ \frac 12 \intQt (|\s|^2+|w|^2)
	\\ & \leq
	\frac \b 2 \I2 {\dt\ph} 
	+ \frac 12 \I2 {\nabla \s}
	+\bigl( \frac {\ch^2} \b + D \hinf +1 \bigr) \I2 \s
	+ \frac {2(1+L^2)}\b \I2 \ph
	\\ & \quad
	+\frac {\ch^2}2 \I2 {\nabla \ph}
	+ \frac 12 \I2 w
	+ \frac {B^2 \ss^2}2 \, T\,|\Omega|.
\end{align*}}
Therefore, a Gronwall argument yields that
\last{
\begin{align}
	\label{First_estimate}
	& \non
	\norma{\m_{\eps}}_{\L\infty H \cap \L2 V}
	+ \norma{\ph_{\eps}}_{\H1 H \cap \L\infty V}
		+ \norma{\s_{\eps}}_{\L\infty H \cap \L2 V}
	\\[0.5mm]
	& \quad
	+ \norma{F_{1,\eps}(\ph)}_{\L\infty {\Lx1}}^{1/2}
	 \, 	\leq\,
	c\,(1+ \norma{u}_{\L2 H}+ \norma{w}_{\L2 H}).
\end{align}%
}%

\noindent {\bf Second estimate:}  \,\,\,
We multiply \eqref{State_second} by $-\Delta \ph$, 
\juerg{write
\last{$F'_\eps=F'_{1,\eps}+F_2'$,}
} 
integrate over $Q_t$, where $t\in(0,T]$, and by parts, to obtain that
\begin{align*}
	&
	\frac \b2 \normaH{\nabla \ph(t)}^2
	+ \I2 {\Delta \ph}
	\last{{}+ \intQt F_{1,\eps}''(\ph)|\nabla \ph|^2}
	\\ & \quad
	=\, 
	\frac \b2 \normaH{\nabla \ph_0}^2
	- \intQt F_2''(\ph)|\nabla \ph|^2
	- \intQt \ch\s\Delta \ph
	- \lastjuerg{\intQt} \m \Delta \ph,
\end{align*}
where the terms on the \rhs~are denoted by $I_1,...,I_4$, in \juerg{this} order.
At first, \juerg{the convexity (recall assumption {\bf (F1)}) of $F_1$ 
\last{and the \wk\ properties of the Yosida regularization 
entail that $ F_{1,\eps}''(\ph)\ge 0$},}
so that the third term on the \lhs~is nonnegative.
Furthermore, the first term $I_1$ on the \rhs\ is bounded 
due to {\bf (A1)}, whereas
the other terms can be dealt \juerg{with by} accounting for Young's inequality and the above estimate.
In fact, we have that
\begin{align*}
	\sum_{i=2}^{4}|I_i| 
	&\, \leq\,
	L \I2 {\nabla \ph}
	+ \ch^2 \I2 {\s}
	+ \I2 {\m}
	+	\frac 12 \I2 {\Delta \ph}.
\end{align*}
Therefore, we realize that $\norma{\Delta \ph}_{\L2 H }^2\leq c$.
The elliptic regularity theory, along with the smooth 
boundary condition in \eqref{State_BC}, and then a comparison in \eqref{State_second},
give us that
\last{
\begin{align}
	\label{second_estimate}
	\norma{\ph_\eps}_{\L2 W}
	+ \norma{F'_ {1,\eps}(\ph_\eps)}_{\L2 H}
	\leq 
	c\,(1+ \norma{u}_{\L2 H}+ \norma{w}_{\L2 H}).
\end{align}
}

\andrea{
\noindent
{\bf Third estimate:} \,\,\,
We now multiply \eqref{State_third} by $\dt\s$, \juerg{and} integrate over $Q_t$ and by parts, to infer that
\begin{align*}
	&\I2 {\dt \s}
	+ \frac B2 \norma{\s(t)}^2_H
	+ \frac 12 \norma{\nabla \s(t)}^2_H
	\\ & \quad 
	= \,
	\frac B2 \norma{\s_0}^2_H
	+ \frac 12 \norma{\nabla \s_0}^2_H
	-\ch \intQt \Delta \ph \,\dt \s
	+ \intQt B\ss\, \dt \s
	- \intQt D \s h(\ph)  \dt \s
	+ \intQt w \,\dt\s.
\end{align*}
Here, it suffices to recall {\bf (A1)}, 
\eqref{First_estimate}, \eqref{second_estimate}, and to
employ Young's inequality several times, to deduce that
\last{
\begin{align}
	\norma{\s_\eps}_{\H1 H \cap \L\infty V}
	\leq 
	c\,(1+ \norma{u}_{\L2 H}+ \norma{w}_{\L2 H}).
\end{align}}
}

\noindent
{\bf Fourth estimate:} \,\,\,
Now, we note that the equation \eqref{State_first} shows a parabolic structure
with respect to $\m$\juerg{; indeed, it can be} rewritten as
\begin{align}
	\begin{cases}
	 \a \dt \m - \Delta \m = f \quad &\hbox{in $Q$}, 
	 \quad \hbox{with} \quad f:= (P\s - A-u)h(\ph) - \dt \ph,
		\\
	\m(0)=\m_0 \quad &\hbox{in $\Omega$}.
	\end{cases} \label{pier5}
\end{align}
On the other hand, owing to the above estimates, the source term \last{$f$ is bounded in} $\L2 H$
and the initial datum is regular, so that the parabolic regularity theory yields that
\begin{align}
	\label{fourth_estimate}
	\norma{\m_\eps}_{\H1 H \cap  \L\infty V \cap \L2 W}
	\leq
	c\,(1+ \norma{u}_{\L2 H}+ \norma{w}_{\L2 H}).
\end{align}

\noindent
{\bf Fifth estimate:} \,\,\,
Next, we differentiate \eqref{State_second} with respect to time and multiply the
resulting equality by $\dt \ph$ to infer that
\begin{align*}
	&
	\frac \b2 \normaH{\dt\ph(t)}^2
	+\I2 {\nabla \dt\ph}
	\last{{}+ \intQt F_{1,\eps}''(\ph)|\dt\ph|^2}
	\\ 
	& \quad
	=
	\frac \b2 \normaH{\dt\ph(0)}^2
	- \intQt F_2''(\ph)|\dt \ph|^2
	+ \intQt \dt\m\,\dt\ph
	+ \intQt \ch\,\dt \s \,\dt \ph.
\end{align*}
\juerg{Again, the third term on the \lhs~is nonnegative.}
On the other hand, the first term on the \rhs~is
under control by virtue of \last{assumptions {\bf (A1)}, {\bf (A2)}, {\bf (F2)} and \eqref{pier4}
which imply  that $F_\eps'(\ph_0)$ is uniformly bounded in $ \Lx\infty$}. In fact,
evaluating \eqref{State_second} at $t=0$, we see that
\begin{align*}
\dt\ph(0)= \frac 1\b [\m_0 + \Delta \ph_0 - \last{F'_\eps(\ph_0)}+\ch\s_0],
\end{align*}
and all \juerg{of the terms on the \rhs\ are bounded} in $\Lx2$.
Lastly, thanks to the Young inequality, we have that
\begin{align*}
	& 
	\sum_{i=2}^{4} |I_i|
	\leq 
	\frac12 \I2 {\dt \m}
	+ \biggl(\frac {1+\ch^2}2 + L \biggr) \I2 {\dt \ph}
	+ \frac12 \I2 {\dt \s}.
\end{align*}
Thus, owing to the previous estimates, we infer that
\last{%
\begin{align}
	\label{fifth_estimate}
	\norma{\ph_\eps}_{W^{1,\infty} (0,T;H) \cap \H1 V }
	\leq
	c\,(1+ \norma{u}_{\L2 H}+ \norma{w}_{\L2 H}).
\end{align}%
}%
\last{Next, by multiplying \last{\eqref{State_second}
by $-\Delta \ph$ and integrating over $\Omega$, we obtain 
\begin{align*}
	 \iO |{\Delta \ph}|^2 
	{}+ \iO F_{1,\eps}''(\ph)|\nabla \ph|^2
	=\, 
    \iO (\b \partial_t \ph  + F_2'(\ph) - \ch \s -  \m )\Delta \ph.
\end{align*}
Then, making use of previous estimates and elliptic regularity results}, 
we easily \last{deduce~that}
\begin{align}
	\label{fifth_estimate_bis}
	\norma{\ph_\eps}_{\L\infty W}
	\leq
	c\,(1+ \norma{u}_{\L2 H}+ \norma{w}_{\L2 H}),
\end{align}}%
which, accounting for the 
Sobolev embedding $\Hx2 \subset \Lx\infty$, also yields that
\last{
\begin{align}
	\label{fifth_estimate_sobolev}
	\norma{\ph_\eps}_{L^\infty (Q)}
	\leq
	c\,(1+ \norma{u}_{\L2 H}+ \norma{w}_{\L2 H}).
\end{align}}

%
%

\noindent
{\bf Sixth estimate:} \,\,\,
\juerg{Next, we observe  that the} equation \eqref{State_third} has parabolic structure with
respect to the variable $\s$, since we can rewrite it as
\begin{align*}
	\begin{cases}
	 \dt \s - \Delta \s + B\s = f \quad &\hbox{in $Q$}, 
	 \quad \hbox{with} \quad f:= \andrea{-\ch \Delta \ph +} B \ss - D \s h(\ph) + w,
		\\
	\s(0)=\s_0 \quad &\hbox{in $\Omega$.}
	\end{cases}
\end{align*}
By virtue of the above estimates and {\bf (A7)}, \last{the reader can easily check} that 
$f$ \last{is bounded in} $\L\infty H$, which allows us to recover 
the full parabolic regularity 
\last{\begin{align}
\label{third_prima}
	\norma{\s_\eps}_{\H1 H \cap \L\infty V \cap \L2 W}
	\leq 
	c\,(1+ \norma{u}_{\L2 H}+ \norma{w}_{\L2 H}).
\end{align}}%
Moreover, provided we assume $\s_0\in \Lx\infty$, as in {\bf (A1)},
we can invoke \cite[Thm.~7.1, p.~181]{LDY} to conclude that
\last{
\begin{align}
	\norma{\s_\eps}_{L^\infty(Q)}
	\leq c\,(1+ \norma{u}_{\L2 H}+ \norma{w}_{\L\infty H}).
	\label{third_estimate}
\end{align}}



\noindent
{\bf Seventh estimate:} \,\,\,
Moreover, the above estimates also entail that the \last{\rhs\ $f$ in \eqref{pier5} 
is bounded in} $\L\infty H$.
By virtue of the assumption $\m_0\in\Lx\infty$, we can again invoke
\cite[Thm.~7.1, p.~181]{LDY} in order to realize that
\last{%
\begin{align}
	\label{sixth_estimate}
	\norma{\m_\eps}_{L^\infty (Q)}
	\leq
	c\,(1+ \norma{u}_{\L\infty H}+ \norma{w}_{\L\infty H}).
\end{align}}%

\noindent
\last{{\bf Passage to the limit:}} \,\,\,
Upon collecting 
\last{the above estimates, all of them independent of $\eps$,
it is a standard matter to realize that 
there is a subsequence of $(\m_\eps, \ph_\eps, \s_\eps)$ suitably 
converging, as $\eps \to 0$, 
to a solution $(\m, \ph, \s)$ of \Stsys\ that verifies
\eqref{stima_reg}. We can argue, for instance, as in \cite{CGH} and just
point out that in the limit procedure we have to use the maximal monotonicity 
of $F_1'$ (intended as a subdifferential operator) and a weak-strong convergence 
argument to identify the limit of $F'_ {1,\eps}(\ph_\eps)$ as  $F'_ {1}(\ph)$.} 

\noindent
{\bf \last{Separation property}:} \,\,\,
At this point, we can rewrite 
the second equation \eqref{State_second} in the form
\begin{align}
	\label{sev_proof}
	\b\dt \ph
	-\Delta \ph
	+ F'(\ph)
	=
	g\last{,}
	 \quad \hbox{with}\quad  g:= \m+\ch \s ,
\end{align}
and, on account of the previous estimates, we \last{know that $g $ is bounded in} $ L^\infty(Q)$,
so that there exists a positive constant $g_*$ for 
which $\norma{g}_{L^\infty(Q)}\leq g_*$.
Besides, the growth assumption {\bf (F4)} implies the existence of some 
constants $r_{*}$ and $r^{*}$ such that $ \rmin < r_{*}\leq r^{*} < \rmax $
and
\begin{gather}
	\label{separation_first}
	 r_{*} \last{{}\leq{}} \infess_{x\in \Omega} \ph_0(x), \quad r^{*} \last{{}\geq{}} \supess_{x\in \Omega} \ph_0(x),
	\\
	F'(r)
	 + g_* \leq 0 
	\quad \forall r \in (\rmin,r_{*}), \quad 
	F'(r)
	 - g_* \geq 0 
	\quad \forall r \in (r^{*},\rmax).
	\label{separation_second}
\end{gather}
Then, let us set, for convenience, $\th:= (\ph-r^*)^+$, multiply equation
\eqref{sev_proof} by $\th$, and integrate over $Q_t$, where $t\in (0,T],$
and by parts, to obtain that
\begin{align*}
	\frac \b2 \normaH{\th(t)}^2 + \I2 {\nabla \th}
	+ \intQt (F'(\ph)-g)\th
	= 0,
\end{align*}
where we also applied \eqref{separation_first} to conclude that $\th(0)=0$.
Moreover, the last term is nonnegative due to 
\eqref{separation_second}, so that $\th = (\ph-r^*)^+ =0$,
which in turn implies that $\ph \leq r^*$ almost everywhere \last{in} $Q$. In a similar
manner, we easily conclude that $\ph \geq r_*$ almost everywhere \last{in} $Q$
by testing \eqref{sev_proof} by $-(\ph-r_*)^-$.
Thus, we have just shown that
\begin{align}
	\label{seventh_estimate}
	r_* \leq {\ph} \leq  r^* \quad \aeQ.
\end{align}
Now, we note that \eqref{separation_result} and 
{\bf (F1)}--{\bf (F4)} directly imply 
\eqref{stima_postsep}.
In fact, \eqref{seventh_estimate} ensures that the phase variable $\ph$
stays away from the boundary of 
\last{$(r_-, r_+)$}, so that
$F$ and its derivatives turn out to be uniformly bounded \last{in $[r_*, r_*]$.}
\Edim

\subsection{Continuous dependence results}
The continuous dependence result to be shown below will
in turn
prove the uniqueness of the solution to the state system \Stsys.

\Bdim[Proof of Theorem~\ref{THM_cont_dep}]
First of all, let us set 
\begin{align}
	\label{diff_variables}
	u:= u_1-u_2, \ w:= w_1-w_2, \ \m:= \m_1-\m_2, \ \ph:= \ph_1-\ph_2, \ \s := \s_1- \s_2.
\end{align}

\last{%
\noindent {\bf First estimate:} \,\,\,
To begin with, we add to both sides of \eqref{State_second} the term $\ph$.
Next, we multiply the difference of \eqref{State_first} by $\m$, 
the difference of the new \eqref{State_second} by $\dt\ph$, and the difference 
of \eqref{State_third} by $\s$. \last{By integrating over $Q_t$, with $t\in (0,T]$,
and adding everything, we obtain a cancellation of terms and arrive at}
\begin{align*} 	
	&
	\frac \a 2 \normaH{ \m(t)}^2
	+ \I2 {\nabla\m}
	+ \b \I2 {\dt\ph}
	+ \frac 12  \normaV{ \ph(t)}^2
	\\ 
	& 
	\quad
	+ \frac 12 \normaH{\s(t)}^2
	+ \intQt B|\s|^2
	+ \intQt |\nabla \s|^2 
	= \,
	-\,\intQt (F'(\ph_1) - F'(\ph_2)) \dt\ph
	\\ 
	& \quad
	+ \intQt \bigl((P\s-u)h(\ph_1)+(P\s_2-A-u_2)(h(\ph_1)-h(\ph_2))\bigr) \m
	\\ 
	& 	\quad
	+ \ch \intQt \s \dt\ph
	+ \intQt  \ph \, \dt\ph 
	\andrea{+\ch \intQt \nabla \ph \cdot \nabla\s}	
	-\intQt D h(\ph_1)|\s|^2
	\\ & \quad
	-\,\intQt D \s_2 (h(\ph_1)-h(\ph_2)) \s
	+ \intQt w\s\,.
\end{align*}
We now estimate the
terms on the \rhs, which we denote by $I_1,...,I_8$, in this order.
We first infer from \eqref{separation_result} that
the nonlinear term $F'$ turns out to be \Lip\ continuous 
\juerg{in the range of interesting arguments},
so that we obtain from Young's inequality that  
\begin{align*}
	- {\intQt (F'(\ph_1) - F'(\ph_2))\dt\ph}
	\leq 
	L \intQt |\ph|  |\dt \ph|
	\leq 
	\frac \b4 \I2 {\dt \ph}
	+ \frac{{L}^2}\b \I2 {\ph},
\end{align*}
where $L$ here stands for a \Lip\ constant of $F'$.
Moreover, accounting for the Young inequality, it is easy to see that
\begin{align*}
	|I_2|
	& \leq
	P \hinf \intQt |\s||\m|	
	+ \hinf \intQt |u||\m|	
	+ \big( P \norma{\s_2}_{L^\infty (Q)} + A + \norma{u_2}_{L^\infty (Q)}\big)L_h \intQt |\ph||\m|	
	\\ &
	\leq
	c \intQt (|\m|^2+|\ph|^2+|\s|^2+|u|^2),
\end{align*}
where we have have used the fact that $\sigma_2$ is a solution to
\Stsys\, and thus has to satisfy \eqref{stima_reg} and also that $u_2$ is an admissible control.
Furthermore, using Young's inequality once more, we have that
\begin{align*}
	|I_3|+|I_4|+|I_5|
	& \leq
	 \frac \b4  \intQt |\dt\ph|^2
	 +\frac 2 \b \intQt(|\s|^2+|\ph|^2)
	+ \frac {\ch^2} 2 \intQt |\nabla \ph|^2
	+ \frac 12 \intQt |\nabla \s|^2.
\end{align*}
Finally, Young's inequality, along with the \Lip\ continuity of $h$, leads \last{us to} 
\begin{align*}
	|I_6|+|I_7|+|I_8|
	& \leq
	D \hinf \intQt |\s|^2
	+	D L_h  \norma{\s_2}_{L^\infty(Q)}\intQt |\ph||\s|
	+ \intQt |w||\s|
	\\ &
	\leq c \intQt (|\ph|^2+|\s|^2+|w|^2).
\end{align*}
\juerg{At this point, we collect the above estimates, and apply  Gronwall's lemma, to conclude that} 
\begin{align}
	\non &
	\norma{\m_1-\m_2}_{\L\infty H \cap \L2 V}
	+ \norma{\ph_1-\ph_2}_{\H1 H \cap \L\infty V}
	\\ & \quad + \norma{\s_1-\s_2}_{\L\infty H \cap \L2 V}
	\leq c \, \bigl(\norma{u_1-u_2}_{\L2 H}+\norma{w_1-w_2}_{\L2 H}\bigr).
	\label{cont_dep_first}
\end{align}}%


\noindent
\last{\bf Second estimate:}
 \,\,\,We multiply the difference of \eqref{State_second} by $-\Delta \ph$,
and use  the Young inequality several times, the previous estimates, and 
the elliptic regularity theory, to obtain that
\begin{align}
	\norma{\ph_1-\ph_2}_{\L2 W}
	\leq c\,\bigl(\norma{u_1-u_2}_{\L2 H}+\norma{w_1-w_2}_{\L2 H}\bigr).
	\label{cont_dep_third}
\end{align}

\noindent
\last{{\bf Third estimate:} }
\,\,\,
Next, we test the difference of \eqref{State_first} by $\dt\m$ and integrate over
time and by parts to realize that
\begin{align*}
	\a \I2 {\dt \m}
	+ \frac 12 \normaH{\nabla \m(t)}^2
	&=
	- \intQt  \dt\ph \dt\m
	+ \intQt (P\s-u)h(\ph_1)\dt\m
	\\ & \quad
	+ \intQt (P\s_2-A-u_2)(h(\ph_1)-h(\ph_2))\dt\m.
\end{align*}
Let us indicate by $I_1,I_2$, and $I_3$ the integrals on the \rhs.
They can be handled, with the help of the Young inequality and the previous estimates, as
follows:
\begin{align*}
	\sum_{i=1}^{3}|I_i| 
	 &\,\leq\, 
	\frac \a2 \I2 {\dt\m}
	+ \frac 1\a\I2 {\dt\ph}
		+ \frac{2 \hinf^2}\a \intQt(P^2 |\s|^2 +|u|^2)
	\\ & \qquad
 	+  L_h \biggl( 
 	P \norma{\s_2}_{L^\infty(Q)} \intQt |{\ph}|\,|{\dt\m}|
 	+ A \intQt |\ph||\dt\m|
 	+ \norma{u_2}_{L^\infty(Q)} \intQt |{\ph}|\,|{\dt\m}|
 	\biggr)
 	\\ &
 	\leq\, 
 	\frac {3\a}4 \I2 {\dt\m}
	+ \frac 1\a \I2 {\dt\ph}
	+ \frac{2 \hinf^2}\a \intQt(P^2 |\s|^2 +|u|^2)
	\\ & \qquad
 	+ \frac{ 3 { L_h}\!^2}\a \biggl( P^2 K_1^2
 	+A^2 + \norma{u_2}^2_{L^\infty(Q)}
 	\biggr)
 	\I2 {\ph},
\end{align*}
where we use the boundedness of $\sigma_2$ once more.
Thus, the above estimates yield that 
\begin{align}
	\label{cont_sep_fourth}
	\norma{\m_1-\m_2}_{\H1 H \cap \L\infty V}
	\leq
	c\,\bigl(\norma{u_1-u_2}_{\L2 H}+\norma{w_1-w_2}_{\L2 H}\bigr).
\end{align}

\noindent
\last{\bf Fourth estimate:} 
\,\,\,
Arguing as in the \last{second} estimate, that is, using  comparison in the 
difference of \eqref{State_first} and elliptic regularity theory, we find that
\begin{align}
	\norma{\m_1-\m_2}_{\L2 W}
	\leq c\,\bigl(\norma{u_1-u_2}_{\L2 H}+\norma{w_1-w_2}_{\L2 H}\bigr).
	\label{cont_dep_fifth}
\end{align}

\noindent
\last{{\bf Fifth estimate:}}
\,\,\,
We multiply the difference of \eqref{State_third} by $\dt\s$, and integrate
over $Q_t$ and by parts, to obtain that
\begin{align*}
	& \I2 {\dt\s}
	+ \frac 12 \normaH{\nabla \s(t)}^2
	+ \frac B2 \normaH{\s(t)}^2
	\\
	 & 	=\, 
	\andrea{-\,\ch \intQt \Delta \ph\, \dt \s}
	- \intQt D \s h(\ph_1)\, \dt \s
	- \intQt D \s_2 (h(\ph_1)-h(\ph_2))\, \dt \s
	+ \intQt w \,\dt \s.
\end{align*}
Here, we denote by \andrea{$I_1,...,I_4$} the terms on the \rhs.
Using Young's inequality four times, along with the \Lip\ continuity
of $h$, we realize that the integrals on the \rhs~can be estimated as follows:
\andrea{
\begin{align*}
 \sum_{i=1}^{4} |I_i|
 & \leq 
	\ch \intQt |\Delta \ph|\,| \dt\s| 	
 	+ D \hinf \intQt |\s| \,|\dt \s|
	\\ & \quad 
	+ D  L_h \norma{\s_2}_{L^\infty(Q)} \intQt |{\ph}|\, |{\dt \s}|
	+ \intQt| w |\,|\dt \s|
 \\ &\leq
 \frac 12 \I2 {\dt \s}
  + 2 \ch^2 \intQt |\Delta\ph|^2
	+ 2D^2 \hinf^2 \I2 {\s}
	\\ & \quad 
 + 2D^2{ L_h}\!^2 K_1^2  \I2 \ph
 + 2\I2 {w},
\end{align*}}
where we again exploit the uniform bound for $\norma{\s_2}_{L^\infty(Q)}$.
Therefore, we deduce that
\begin{align}
	\label{cont_dep_sixth}
	\norma{\s_1-\s_2}_{\H1 H \cap \L\infty V}
	\leq c\,\bigl(\norma{u_1-u_2}_{\L2 H}+\norma{w_1-w_2}_{\L2 H}\bigr).
\end{align}

\noindent
\last{{\bf Sixth estimate:} } \,\,\,
Finally, by comparison in the difference of \eqref{State_second}, and 
applying elliptic regularity theory, we have that
\begin{align}
	\norma{\s_1-\s_2}_{\L2 W}
	\leq c\,\bigl(\norma{u_1-u_2}_{\L2 H}+\norma{w_1-w_2}_{\L2 H}\bigr).
	\label{cont_dep_seventh}
\end{align}
Upon collecting all of the estimates \accorpa{cont_dep_first}{cont_dep_seventh}, 
we find that \eqref{cont_dep_final}
is shown, so that Theorem~\ref{THM_cont_dep} is completely proved.
\Edim

\section{The control problem}
\label{SEC_CONTROL_PROBLEM}
\setcounter{equation}{0}
From now on, we are going to focus our attention
on the control problem. The main results are the existence
of optimal controls and the first-order necessary
conditions for optimality.

\subsection{Existence of optimal controls}
\Bdim[Proof of Theorem~\ref{THM_existence_opt_control}]
The proof makes use of the direct method from the calculus of variations.
In fact, the cost functional is nonnegative, convex, and weakly lower semicontinuous.
\juerg{To this end, let us pick a minimizing sequence
$\graffe{(u_n,w_n)}_{n\in\enne}\subset \Uad$ such that, setting
$(\m_n,\ph_n,\s_n)=\S(u_n,w_n)$, and recalling the notations \accorpa{cost}{Jred}, }
there holds
\begin{align*}
	\lim_{n\to \infty} \J (\ph_n,\s_n,u_n,w_n) 
	=\lim_{n\to \infty} \Jred (u_n,w_n) 
	= \inf_{(u,w)\in \Uad} \Jred (u,w) .
\end{align*}
On the other hand, $\graffe{(u_n,w_n)}_{n\in\enne}$ is bounded in $L^\infty(Q)\times L^\infty(Q)$, 
and also the bounds \eqref{stima_reg} and \eqref{stima_postsep} 
are at our disposal\last{, for every $n\in\enne$ \last{and for the 
corresponding triplet} $(\m_n,\ph_n,\s_n)$}.
Hence, accounting for standard weak compactness arguments 
(see, e.g., \cite[Sec.~8, Cor.~4]{Simon}), it is a standard matter to
infer the existence of a pair $(\uopt,\wopt)$ and a triplet $(\bm,\bph,\bs)$ such that
the following convergence properties are (possibly only on a subsequence) fulfilled as $n \to \infty$:
\begin{align}
	(u_n,w_n) & \to (\uopt,\wopt) \quad \hbox{weakly star in} \ (L^\infty (Q))^2,
	\\
	\non
	\m_n & \to \bm \quad \hbox{weakly star in} \ {\H1 H \cap \L\infty V \cap \L2 W \cap L^\infty (Q)}
	\\ 
	& \hspace{2.2cm} \hbox{and strongly in} \ \C0 H \cap \L2 V,
	\\\non
	\ph_n & \to \bph \quad \hbox{weakly star in} \ W^{1,\infty}(0,T;H) \cap \H1 V \cap \L\infty W
	\\ 
	& \hspace{2.2cm} \hbox{and strongly in} \ C^0(\overline Q),
	\label{phntobph}
	\\ \non
	\s_n & \to \bs \quad \hbox{weakly star in} \ {\H1 H \cap \L\infty V \cap \L2 W \cap L^\infty (Q)}
	\\	& \hspace{2.2cm} \hbox{and strongly in} \ \C0 H \cap \L2 V.
\end{align}
Clearly, as the convex set $\Uad$ is weakly sequentially closed, we have that $(\uopt,\wopt)\in \Uad$;
besides, the \juerg{strong convergence properties show that the Cauchy conditions \eqref{State_IC} are fulfilled by
$(\bm,\bph,\bs)$.} 
Moreover, the strong convergence in \eqref{phntobph} and the 
assumptions {\bf (F1)}--{\bf (F4)} and {\bf (A3)} imply that
\begin{align*}
	h(\ph_n) \to h(\bph) \ \  \hbox{and} \ \  F'(\ph_n) 
	\to F'(\bph) \quad  \hbox{strongly in} \ C^0(\overline Q),\quad\mbox{as \,$n\to\infty$}.
\end{align*}
\juerg{Therefore, passing to the limit as $n\to\infty$ in the corresponding time-integrated version
of \Stsys, written for $(u_n,w_n)$ and $(\m_n,\ph_n,\s_n)$, we easily see that}
$(\bm,\bph,\bs)$ solves \Stsys\ with $(u,w)=(\uopt,\wopt)$, which
yields that $(\bm,\bph,\bs)=\S(\uopt,\wopt)$.
Finally, we combine the weak sequential lower semicontinuity of the cost functional
with the assumption that $(u_n,w_n)$ is a minimizing sequence 
to deduce that $(\uopt,\wopt)$ is indeed an optimal control.
\Edim

\subsection{The linearized system}
%
At this point, our aim is to find the
necessary conditions for optimality. Actually, we would like to
\juerg{express} the formal variational inequality \eqref{formal_foc} in an explicit form.
\juerg{For this purpose}, we have to prove the \Frechet\ differentiability of the
reduced cost functional $\Jred$, which is the 
composition of $\J$ with the control-to-state operator $\S$.
However, $\J$ is \sfw ly \Frechet\ differentiable. Therefore, it suffices to 
prove that $\S$ is \Frechet\ differentiable as well, and then invoke the chain rule to 
write \eqref{formal_foc} in an explicit way.

The expectation is that, provided we find the proper Banach spaces,
the \Frechet\ derivative of $\S$ applied to the pair $(k,l)$
is given by the unique solution
to the linearized system \Linsys. With this in mind, we begin 
by establishing the well-posedness of the linearized system \Linsys.

\Bdim[Proof of Theorem~\ref{THM_well-posedness_lin}]
For the sake of simplicity, we proceed formally
\last{by providing some estimates. 
Here, let us just mention that a 
rigorous proof can be carried out, e.g.,
using a Galerkin scheme: first showing the rigorous counterpart of the following
estimates so to ensure that they are
independent of the discretization parameter, and then passing to the limit 
with respect to the discretization parameter.}
Moreover, the system \Linsys\ is linear,
so that the uniqueness directly follows from the uniform estimates.
In addition, some of the forthcoming estimates follow the same lines 
as the ones of the state system, which \juerg{allows} us to be less detailed
below.

\andrea{
\noindent
{\bf First estimate:} \,\,\, 
First of all, we add to both sides of \eqref{Lin_second} the term $\xi$.
Then, we multiply \eqref{Lin_first} by $\et$, the new \eqref{Lin_second} by $\dt\xi$,
\eqref{Lin_third} by $\z$, \juerg{add the resulting equations,} and integrate over $Q_t$ and by parts
for an arbitrary $t\in(0,T]$. After a cancellation \juerg{of terms} and some rearrangements, and making use of
the initial conditions 
\eqref{Lin_IC}, we obtain that
\begin{align*}
	& 
	\frac \a2 \normaH{\et(t)}^2
	+ \I2 {\nabla \et}
	+ \b \I2 {\dt \xi}
	+ \frac 12 \normaV{\xi(t)}^2
	+ \frac 12 \normaH{\z(t)}^2
	\\ &\quad
	+ B \I2 \z
	+ \I2 {\nabla \z}
	=
	\intQt (P\z - k) h(\bph)\et
	+ \intQt ( P\bs - A -\uopt) h'(\bph)\xi \et
	\\ &\quad
	+ \intQt \xi \, \dt\xi
	- \intQt F''(\bph)\xi \dt\xi
	+ \ch \intQt  \z\dt\xi
	\andrea{+\ch\intQt \nabla \xi \cdot \nabla \z}
	- \intQt D h(\bph) |\z|^2
	\\ &\quad
	- \intQt D \bs h'(\bph)\xi \z
	+\intQt  l \z\,.
\end{align*}
\juerg{We denote  by $I_1,...,I_9$ the} integrals on the \rhs.
Using the Young inequality, we infer that
\begin{align*}
	|I_1|+|I_2|
	& \leq
	P\hinf \intQt  |\z||\et|
	+ \hinf \intQt  |k||\et|
	+ h'_\infty  \bigl( P \norma{\bs}_{L^\infty(Q)} + A + \norma{\uopt}_{L^\infty(Q)}\bigr)
		\intQt |\xi||\et|
	\\ &
	\leq
	\frac {P\hinf }2 \intQt  (|\z|^2+|\et|^2)
	+\frac {\hinf }2 \intQt  (|k|^2+|\et|^2)
	\\ & \quad
	+ \frac {h'_\infty  \bigl( P K_1 + A + \norma{u^*}_{L^\infty(Q)}\bigr)}2
	\intQt  (|\xi|^2+|\et|^2),
\end{align*}
where we use the fact that $\bs$ satisfies \eqref{stima_reg} and $\uopt$ belongs to 
the class of admissible controls.
Moreover, from Young's inequality, combined with \eqref{stima_postsep},
it follows that
\begin{align*}
	\sum_{i=3}^{5}|I_i| 
	& \leq
	\intQt |\xi||\dt\xi|
	+\norma{F''(\bph)}_{L^\infty(Q)}	\intQt |\xi||\dt\xi|
	+\ch \intQt |\z||\dt\xi|
	\\ & 
	\leq
	\frac \b 2 \I2 {\dt\xi}
	+ \frac 3{2\b} \bigl( K_2^2 +1 \bigr) \I2 {\xi}
	+ \frac {3\ch^2}{2\b}  \I2 {\z}	
\end{align*}
and also that
\begin{align*}
	\sum_{i=7}^{9}|I_i| 
	& \leq
	D \hinf \intQt  |\z |^2
	+ Dh'_\infty \norma{\bs}_{L^\infty(Q)} \intQt |\xi||\z|
	+\intQt  |l||\z|
	\\ & 
	\leq
	\biggl( D\hinf + \frac{(D \hinf' K_1)^2+1}4 \biggr) \I2 \z
	+ \intQt (|\xi|^2+| l|^2).
\end{align*}
Furthermore, using Young's inequality once more, we infer that
\begin{align*}
	|I_6|
	\leq 
	\frac 12 \I2 {\nabla \z}
	+ \frac {\ch^2}2 \I2 {\nabla \xi}.
\end{align*}
\juerg{At this point, we collect all of the above estimates} and apply Gronwall's lemma to deduce that
\begin{align}
	\non 
	& \norma{\et}_{\L\infty H \cap \L2 V}
	+ \norma{\xi}_{\H1 H \cap \L\infty V}
	+ \norma{\z}_{\L\infty H \cap \L2 V}
	\\[0.5mm]
	& 	\leq\, c\,(\norma{k}_{L^2(Q)}+\norma{l}_{L^2(Q)}).
	\label{linearized_first}
\end{align}	
}

\noindent
{\bf Second estimate:} \,\,\,
\juerg{We now observe that the equation} \eqref{Lin_first} shows a parabolic
structure with respect to the variable $\et$. In fact, we can write 
\eqref{Lin_first} in the form 
\begin{align*}
	\a \dt \et
	- \Delta \et
	=f_1 \quad \hbox{with}\quad 
	f_1:=	(P\z - k) h(\bph)
	+( P\bs - A -\uopt) h'(\bph)\xi
	-  \dt \xi,
\end{align*}
\juerg{where}, owing to the above estimate, we easily \juerg{verify} that $f_1\in\L2 H$
and 
\begin{align}
	\label{pier1}
	\norma{f_1}_{\L2 H}
	\,\leq\, 
	c\,(\norma{k}_{L^2(Q)}+\norma{l}_{L^2(Q)}).
\end{align}
So, recalling \juerg{the} boundary and initial conditions \accorpa{Lin_BC}{Lin_IC},
it is a standard matter to recover the full parabolic regularity
and infer that
\begin{align}
	\label{linearized_second}
	\norma{\et}_{\H1 H \cap \L\infty V \cap \L2 W}
	\,\leq\,
	c\,(\norma{k}_{L^2(Q)}+\norma{l}_{L^2(Q)}).
\end{align}

\andrea{
\noindent
{\bf Third estimate:} \,\,\,
In the same way, we also have
\begin{align*}
	 \b\dt \xi
	-\Delta \xi
	&= f_2 \quad \hbox{with}\quad
	f_2:= - F''(\bph)\xi
	+ \ch \z
	+\et,
	\\[1mm]
	\dt \z
	- \Delta \z
	&=f_3
	\quad \hbox{with}\quad
	f_3:=
	\andrea{-\ch\Delta \xi}
	-B\z
	- D \z h(\bph)
	- D \bs h'(\bph)\xi
	+ l.
\end{align*}
Then, we first note that $f_2$ belongs to $\L2 H$ and satisfies the same estimate
as in \eqref{pier1}, so that 
the regularity theory for parabolic equation with regular initial datum and
homogeneous Neumann boundary conditions allows us to to infer that
\begin{align}
	\label{linearized_third}
	\norma{\xi}_{\H1 H \cap \L\infty V \cap \L2 W}
	\,\leq\,
	c\,(\norma{k}_{L^2(Q)}+\norma{l}_{L^2(Q)}).
\end{align}
\juerg{Besides}, also $f_3$ belongs to $\L2 H$, and similar reasoning leads to \juerg{the conclusion} that
\begin{align}
	\label{linearized_third_bis}
	\norma{\z}_{\H1 H \cap \L\infty V \cap \L2 W}
	\,\leq\,
	c\,(\norma{k}_{L^2(Q)}+\norma{l}_{L^2(Q)}),
\end{align}
which concludes the proof of Theorem~\ref{THM_well-posedness_lin}.}
\Edim

\subsection{Differentiability of the control-to-state operator}
Now we are going to show the \Frechet\ differentiability of the operator
$\S$ and \juerg{to} characterize its \Frechet\ derivative.
\Bdim[Proof of Theorem~\ref{THM_frechet}]
At first, let us fix a control pair $(\uopt,\wopt)\in \Uad\subset \UR$ with the corresponding state
$(\bm,\bph,\bs)$. Then, whenever $(k,l)$ belongs to $(L^\infty(Q))^2$, we denote with $(\et,\xi,\z)$
the corresponding solution to system \Linsys. 
Moreover, let us recall that $\UR$ is an open set, so that, provided
that we consider small perturbations, we also have $(\uopt+k, \wopt+l)\in\UR$.
Namely, there exist some positive constant $\d_*$ such that
$(\uopt+k, \wopt+l)\in\UR$ for every $(k,l)$ such that 
$\norma{k}_{L^\infty(Q)}+\norma{l}_{L^\infty(Q)}\leq \d_*$.
In the following, we \juerg{always} assume that this is the case.
Lastly, we denote with $(\hm,\hph,\hs)$ the unique solution to \Stsys\ corresponding
to the incremented control $(\uopt+k, \wopt+l)$.
Let us point out that Theorem~\ref{THM_well-posedness_lin}
entails that the map $(k,l)\mapsto(\et,\xi,\z)$ is linear and continuous
between $(L^2(Q))^2$ and $(\H1 H \cap \L\infty V \cap \L2 W)^3$.

Here, we aim \juerg{at directly checking} the definition of \Frechet\ differentiability
for $\S$. Namely, we \juerg{are going to show} that
\begin{align}
	&\non
	\S (\uopt+k,\wopt+l) = \S (\uopt,\wopt) + [D{\S}(\uopt,\wopt)](k,l) + o(\norma{(k,l)}_{L^\infty(Q)\times L^\infty(Q)})
   \quad \hbox{in $\cal Y$} \ \  \\[1mm]
    & 
      \hbox{as} \ \  \norma{(k,l)}_{L^\infty(Q)\times L^\infty(Q)}\to 0,
   \label{fre_formal_one}
\end{align}
for the Banach space  \,$\cal Y$\, \juerg{introduced} in \eqref{pier2}.
\juerg{To this end}, it is convenient to set
\begin{align*}
	\ps:= \hm - \bm-\et, \ \ y:= \hph-\bph-\xi,\  \ z:=\hs-\bs-\z.
\end{align*}
\juerg{With this notation, \eqref{fre_formal_one} takes the form}
\begin{align*}
	&\norma{(\ps,y,z)}_{\cal Y} = o(\norma{(k,l)}_{L^\infty(Q)\times L^\infty(Q)})
	\quad
   \hbox{as} \ \  \norma{(k,l)}_{L^\infty(Q)\times L^\infty(Q)}\to 0.
\end{align*}
\juerg{Obviously, the validity of this condition implies that $\,\S\,$ is \Frechet\ differentiable at $\,(\uopt,\wopt)\,$ and that
$[D{\S}(\uopt,\wopt)](k,l)=(\et,\xi,\z)$ for every $(k,l)\in (L^\infty(Q))^2$.
To verify this condition, it suffices to construct an increasing function $\,G:(0,\d_*)\to (0,+\infty)\,$ such that
$\norma{(\ps,y,z)}_{\cal Y}^2 \leq G \bigl(\norma{(k,l)}_{L^\infty(Q)\times L^\infty(Q)}\bigr)$ and
\begin{align}
	\label{fre_formal}
	\lim_{\lam \to 0} \frac {G(\lam)}{\lam^2}= 0.
\end{align}
This is actually the estimate we are going to check with the
choice $G(\lam)=c \lam^4$ for some positive constant $c$.}

At this stage, let us recall that since $(\hm,\hph,\hs)$ and $(\bm,\bph,\bs)$ are fixed,
they both verify \eqref{stima_reg} and \eqref{stima_postsep},
as well as the following continuous dependence estimate
\begin{align}
	\non &
	\norma{\hm-\bm}_{\H1 H \cap \L\infty V \cap \L2 W}
	+ \norma{\hph-\bph}_{\H1 H \cap \L\infty V \cap \L2 W}
	\\[0.5mm]
	 &   \qquad \non
	+ \norma{\hs-\bs}_{\H1 H \cap \L\infty V \cap \L2 W}
	\\[0.5mm] 
	& \quad
	\leq\, 
	K_3 \, \bigl(\norma{k}_{\L2 H}+\norma{l}_{\L2 H}\bigr),
	\label{fre_contdep}
\end{align}
which directly follows from \eqref{cont_dep_final}.

Besides, a system for $(\ps,y,z)$ can be constructed in light of the
systems \Stsys~corresponding to $(u,w)= (\uopt+k,\wopt+l)$, 
\Stsys~for $(u,w)= (\uopt,\wopt)$,
and \Linsys. By combining them, we obtain the following system\lastjuerg{:}
\begin{align}
	\label{Fre_first}
	& \non
	\a \dt \ps
	+ \dt y
	- \Delta \ps
	=
	Pz h(\bph)
	+ (P\bs -A - \uopt)(h(\hph)-h(\bph)-h'(\bph)\xi)
	\\[0.5mm] & \hspace{3.5cm}
	-k (h(\hph)-h(\bph)) 
	+ P(\hs-\bs)(h(\hph)-h(\bph))
	\quad\hbox{\,in $\,\, Q$},
	\\[1mm]
	\label{Fre_second}
	&\ps
	=
	\b\dt y
	-\Delta y
	+ (F'(\hph) - F'(\bph) - F''(\bph)\xi)
	- \ch z
 	\quad\hbox{\,in $\,Q$},
	\\[1mm]
	\non
	&\dt z
	- \Delta z
	+Bz
	=
	\andrea{-\ch \Delta y}
		- D [ \bs (h(\hph)-h(\bph)-h'(\bph)\xi)
	\\[1mm]
	 &
	\hspace{3.2cm}
	+ (\hs-\bs)(h(\hph)-h(\bph))
	+ h(\bph)z ]
	\quad \hbox{\,in $\,\,Q$},
	\label{Fre_third}
	\\[1mm]
	&\dn \ps=\dn y=\dn z =0
 	\quad \hbox{\,on $\,\,\Sigma$},
	\label{Fre_BC}
	\\[1mm]
	& \ps(0)=y(0)=z(0)=0
  	\quad\hbox{\,in $\,\Omega$}.
	\label{Fre_IC}
\end{align}
\Accorpa\Fresys Fre_first Fre_IC

Note that \regstate~and \eqref{reg_lin} entail that
\begin{align*}
	\ps,y,z \in \H1 H \cap \L\infty V \cap \L2 W.
\end{align*}

\noindent
{\bf First estimate:} \,\,\,
First of all, we add to both sides of \eqref{Fre_second} 
the term $y$.
Next, we multiply \eqref{Fre_first} by $\ps$, 
the new \eqref{Fre_second} by $\dt y$,
and \eqref{Fre_third} by $z$. Then, we add \juerg{the resulting identities}, integrate over $Q_t$,
where $t\in(0,T]$, and by parts, to find  that
\begin{align*}
	&
	\frac \a2 \normaH{\ps(t)}^2
	+ \I2 {\nabla \ps}
	+ \b \I2 {\dt y}
	+ \frac 12 \normaV{y(t)}^2
	+ \frac 12 \normaH{z(t)}^2
	\\ & \qquad
	+ B\I2 {z}
	+ \I2 {\nabla z}
	\\ & \quad
	=
	\intQt P z h(\bph) \ps
	+ \intQt(P\bs -A - \uopt)(h(\hph)-h(\bph)-h'(\bph)\xi) \ps
	\\ & \qquad
	-\intQt k (h(\hph)-h(\bph))  \ps
	+ \intQt P(\hs-\bs)(h(\hph)-h(\bph)) \ps
	\\ & \qquad
	- \intQt (F'(\hph) - F'(\bph) - F''(\bph)\xi) \dt y
	+ \intQt \ch z \dt y
	\\ & \qquad
	+ \intQt y \dt y
	\andrea{+\ch\intQt \nabla y\cdot \nabla z}
	-\intQt D \bs (h(\hph)-h(\bph)-h'(\bph)\xi)z
	\\ & \qquad
	-\intQt D (\hs-\bs)(h(\hph)-h(\bph))z
	-\intQt D h(\bph)|z|^2,
\end{align*}
where we denote \andrea{by $I_1,...,I_{11}$} the integrals on the \rhs. Moreover, in 
the above calculations we also owe to the fact that the initial data are zero by \eqref{Fre_IC}.
Using the \Holder\ and Young inequalities, the \Lip\ continuity of $h$
and the Sobolev embedding \eqref{VinL6} with $q=4$, we have that
\begin{align*}
	|I_1|+|I_3|+|I_4|
	&
	\leq
	P\hinf \intQt |z||\ps|
	+ L_h \iot \norma{k(s)}_2\norma{\hph(s)-\bph(s)}_4\norma{\ps(s)}_4\,ds
	\\ & \quad
	+ P L_h \iot \norma{\hs(s)-\bs(s)}_4\norma{\hph(s)-\bph(s)}_4\norma{\ps}_2\,ds
	\\ & \leq 
	\frac {P\hinf}2 \intQt (|{z}|^2+|\ps|^2)
	+ \frac 12 \iot \norma{\ps(s)}^2_V\,ds
	+ c\, \norma{\hph-\bph}_{\L\infty V}^2 \I2 {k}
	\\ & \quad
	+ c \, \norma{\hs-\bs}_{\L\infty  {V}}^2 \norma{\hph-\bph}_{\L\infty {V}}^2 
	+ \last{\I2 {\ps}}
	\\ & 
	\juerg{\leq 
	\frac 12 \iot \norma{\ps(s)}^2_V\,ds
	+ c \,\bigl(\norma{k}_{\L2 H}^4+\|l\|_{L^2(0,T;H)}^4\bigr)
	{+ \,c \intQt (|z|^2 + |\ps|^2 )}}\,,
\end{align*}
where we also invoked the continuous dependence estimate \eqref{fre_contdep}.
Before moving on, let us recall the Taylor formula with integral remainder
which will be useful to estimate some terms.
For an arbitrary function $g\in C^1(\erre)$ with 
$g'$ \Lip\ continuous, we have that
\juerg{
\begin{align}
	\label{taylor_formal_after}
	g(x)= g(\overline x) + g'(\overline x) (x - \overline x) 
	+  (x - \overline x)^2  \int_{0}^1 g''( \overline x + s (x - \overline x))(1-s) \, ds
	\quad \hbox{ for every $x \in \erre$}.
\end{align}
}
Applying the above formula to $F'$ and $h$, respectively, we infer that
\begin{align}
	\label{taylor_F}
	F'(\hph) - F'(\bph) - F''(\bph)\xi \, &=\, F''(\bph)y + R_1 (\hph-\bph)^2,
	\\ 
	\label{taylor_h}
	h(\hph)-h(\bph)-h'(\bph)\xi \, &=\, h'(\bph) y  + R_2 (\hph-\bph)^2,
\end{align}
\juerg{
with the remainders
\begin{align*}
	R_1:= \int_0^1 F'''(\bph+s (\hph-\bph)) (1-s)ds\,, 
	\quad
	R_2:= \int_0^1 h''(\bph+s (\hph-\bph)) (1-s)ds\,.
\end{align*}
Taking advantage of \eqref{stima_postsep} and {\bf (A3)}, we see that
\begin{align*}
	\norma{R_1}_{L^\infty(Q)} \leq R_1^*, \quad
	\norma{R_2}_{L^\infty(Q)}\leq R_2^*,
\end{align*}
for some positive constants $R_1^*,R_2^*$.}
Thus, making use of \eqref{taylor_h}, we are now in a position
to estimate $I_2$ as follows:
\begin{align*}
 	|I_2| 
 	&\leq
 	(P\norma{\bs}_{L^\infty(Q)} +A + \norma{\uopt}_{L^\infty(Q)})
 	\intQt(\hinf'|y| + R_2^* (\hph-\bph)^2) |\ps|
 	\\ & 
 	\leq 
 	\frac {(P K_1+A+\norma{\uopt}_{L^\infty(Q)})h'_\infty }{2} 
 	\intQt (|y|^2+|\ps|^2)
  	\\ & \quad
 	+ (P K_1+A+\norma{\uopt}_{L^\infty(Q)})R_2^* 
 	\iot \norma{\hph(s)-\bph(s)}_4^2\, \norma{\ps(s)}_2\,ds
   	\\ & \leq
	\frac {(P K_1+A+\norma{\uopt}_{L^\infty(Q)})h'_\infty }{2} 
 	\intQt (|y|^2+|\ps|^2)
 	\\ & \quad
 	+ c \norma{\hph-\bph}_{\L\infty {V}}^4
 	+ \I2 \ps
 	\\ &  
 	\leq
 	c \intQt (|y|^2+|\ps|^2) + c (\norma{k}^4_{L^2(Q)}+\norma{l}^4_{L^2(Q)})\,,
\end{align*}
where we also use \eqref{fre_contdep}, the fact that $\bs$ is bounded for \eqref{stima_reg},
whereas $\uopt$ is bounded since it is an admissible control.
As for $I_5$, thanks to the Young inequality and \eqref{taylor_F}, we have that
\begin{align*}
	|I_5|
	&\leq
	\frac \b4 \I2 {\dt y}
	+ \frac{2\norma{F''(\bph)}_{L^\infty(Q)}^2}\b \I2 {y}
	+ 2{R_1^*}^2 \norma{\hph-\bph}_{\L\infty {\Lx4}}^4
	\\ & 
	\leq 
	\frac \b4 \I2 {\dt y}
	+ \frac {2 K_2^2}\b \I2 {y}
	+ c\, (\norma{k}^4_{L^2(Q)}+\norma{l}^4_{L^2(Q)}).
\end{align*}
Moreover, using the Young inequality once more, we have that
\begin{align*}
	\sum_{i=6}^{8} |I_i|	
	\leq
	\frac 12 \I2 {\nabla z}
	+ \frac {\ch^2}2 \I2 {\nabla y}
	+ \frac \b4 \I2 {\dt y}
	+ \frac {2\ch^2}{\b} \I2 {z}
	+ \frac 2{\b} \I2 {y}.
\end{align*}
Lastly, \juerg{by} similar reasoning, we obtain that
\begin{align*}
	\sum_{i=9}^{11} |I_i|	
	& \leq
	D \norma{\bs}_{L^\infty(Q)} \intQt (h'_\infty |y|  + R_2^* (\hph-\bph)^2)|z|
	\\ & \quad
	+ D  L_h \iot \norma{\hs(s)-\bs(s)}_4\,\norma{\hph(s)-\bph(s)}_4\,\norma{z(s)}_2\,ds
	+ D \hinf \I2 {z}
	\\ &
	\leq
	\frac {D K_1 h'_\infty}2 \intQt (|y|^2+|z|^2)
	+D K_1 R_2^* \iot \norma{\hph(s)-\bph(s)}_4^2\, \norma{z(s)}_2\,ds
	\\ & \quad
	+ \frac{D^2  L_h^2}4 \norma{\hs-\bs}_{\L\infty {\Lx4}}^2 \norma{\hph-\bph}_{\L\infty {\Lx4}}^2 
	+\I2 {z}
	\\ & \quad
	+ D \hinf \I2 {z}
	\\ &
	\leq
	\frac {D K_1 h'_\infty}2 \intQt (|y|^2+|z|^2)
	+ c \norma{\hph-\bph}_{\L\infty {V}}^4 
	\\ & \quad
	+ c \norma{\hs-\bs}_{\L\infty {V}}^2\norma{\hph-\bph}_{\L\infty {V}}^2
	+ (D \hinf +1) \intQt  |z|^2
	\\ &
	\leq
	c \intQt (|y|^2+|z|^2)
	+ c (\norma{k}^4_{L^2(Q)}+\norma{l}^4_{L^2(Q)}).
\end{align*}
Hence, applying Gronwall's lemma, we deduce that
\begin{align*}
	&
	\norma{\ps}^2_{\C0 H \cap \L2 V}
	+\norma{y}^2_{\H1 H \cap \L\infty V }
	+\norma{z}^2_{\C0 H \cap \L2 V}
	\\[1mm]
	 & 
	\leq\, 
	C
	\norma{(k,l)}^4_{L^2(Q)\times L^2(Q)},
\end{align*}
which in turn implies \juerg{the validity of \eqref{fre_formal} 
with the \lastjuerg{choice} $G(\lam)=C \lam^4$.
This concludes the proof of the assertion as it suffices to employ the continuous embedding 
$L^\infty(Q)\times L^\infty(Q) \subset L^2(Q)\times L^2(Q)$ to conclude.}

\Edim

\subsection{First-order necessary optimality conditions}

As already pointed out in Section~\ref{SEC_MATHEMATICAL_SETTING}, 
we would like to employ the adjoint variables
in order to eliminate the linearized variables from the variational inequality \eqref{foc_first}.
Here, we begin with the task of establishing the well-posedness of the adjoint system.
In this direction, let us set
\begin{align*}
	Q_t^T \:= (t,T)\times \Omega \quad \hbox{for every $t\in(0,T)$}.
\end{align*}

\Bdim[Proof of Theorem~\ref{THM_wellposed_adjoint}]
The rigorous proof should employ an approximation technique. 
Anyhow, since the system is linear and the arguments are
standard, we simply point out the estimates which allow us to conclude,
leaving the details to the reader. It is worth recalling that the adjoint 
system is linear, so that the uniqueness directly follows from our estimates.

{\bf First estimate:}
First, we add to both sides of \eqref{Adj_first} the term $-q$.
Then, we multiply the new \eqref{Adj_first} by $-\dt q$, \eqref{Adj_second} by $p$,
\eqref{Adj_third} by $\ch^2 r$, add the resulting equations, and integrate over $Q_t^T$ and by parts.
We obtain a cancellation and deduce that
\begin{align*}
	&
	\a \intQtT |\dt q |^2
	+ \frac 12 \normaV{q (t)}^2
	+ \frac \b2  \normaH{p(t)}^2
	+ \intQtT |\nabla p|^2
	+ \frac {\ch^2 }2  \normaH{r(t)}^2
	\\ & \qquad 
	+ \ch^2 \intQtT |\nabla r|^2
	+\ch^2 B \intQtT |r|^2
	\\ & \quad 
	=
	\frac {\bO^2}{2\b}  \normaH{\bph(T)-\phO}^2
	+\frac {{\ch}^2 {\bG}\!^2}2 \normaH{\bs(T)-\sO}^2
	+ \intQtT q \dt q
	\\ & \qquad
	\andrea{+\ch \intQtT \nabla r \cdot \nabla p}
	- \intQtT F''(\bph)|p|^2
	+ \intQtT (P\bs -A -\uopt)h'(\bph) q p
	\\ & \qquad
	- \intQtT D\bs h'(\bph) r p
	+\intQtT \bQ (\bph-\phQ) p
	- \ch^2 \intQtT D h(\bph)|r|^2
	\\ & \qquad
	+ \ch^3 \intQtT p r
	+ \ch^2  \intQtT P h(\bph)q r
	+\ch^2  \intQtT  \bS (\bs - \sQ) r,
\end{align*}
where we used the information \eqref{Adj_IC} on the final data.
In the above equality, the terms on the \lhs~are nonnegative, whereas 
we denote the integrals on the \rhs~by \andrea{$I_1,...,I_{12}$},
in this order. 
As far as the \rhs~is concerned, the 
first four terms can be easily handled with the aid of  \eqref{stima_reg},
assumption \andrealast{{\bf (A6)}}, and the 
Young inequality. Indeed, we have
\begin{align*}
	\sum_{i=1}^{4}|I_i| 
	\leq
	c + \frac \a2 \intQtT |\dt q|^2
	+ \frac 1 {2\a} \intQtT |q|^2
	+ \frac 1 {2} \intQtT |\nabla p|^2
	+ \frac {\ch^2} {2} \intQtT |\nabla r|^2.
\end{align*}
Using Young's  inequality, we can deal with $I_6$ as follows:
\begin{align*}
	|I_6|
	\leq
	\frac{ ( P \norma{\bs}_{L^\infty (Q)} + A 
	+\norma{\uopt}_{L^\infty (Q)} )h'_\infty}2 \intQtT (|q|^2+ |p|^2),
\end{align*}
where we employ that $\bs$ satisfies \eqref{stima_reg} and
that $\uopt$ is an admissible control.
The rest of the terms can be handled using several times the Young inequality
to get that
\begin{align*}
	|I_5|
	+ \sum_{i=7}^{12}|I_i| 
	& \leq
	\frac {P\hinf{\ch^2} }2 \intQtT |q|^2
	+ \biggl( \norma{F''(\bph)}_{L^\infty (Q)} 
	+ \frac {2 + D\norma{\bs}_{L^\infty (Q)} h'_\infty +\ch^3 }2\biggr) \intQtT |p|^2
	\\ & \qquad
	+ \biggl( \frac {2 + D\norma{\bs}_{L^\infty (Q)} h'_\infty +\ch^3 }2 
	+ {\ch^2}D\hinf + \frac {P\hinf{\ch^2}}2 \biggr) \intQtT |r|^2
	\\ & \qquad
	+ \frac {(\bQ^2 + {\ch^4} \bS^2)}4 \intQtT (|\bph - \phQ|^2+ |\bs - \sQ|^2).
\end{align*}
Thus, the backward-in-time Gronwall lemma yields that
\begin{align}
\label{adjoint_first}
	& 
	\norma{q}_{ \H1 H \cap \L\infty V }
	+ \norma{p}_{\L\infty H \cap \L2 V}
	+ \norma{r}_{\L\infty H \cap \L2 V}
	\leq c.
\end{align}

\noindent
{\bf Second estimate:} \,\,\,
By \eqref{adjoint_first} and a comparison argument in \eqref{Adj_second} and \eqref{Adj_third}, 
we obtain that
\begin{align}
\label{adjoint_second}
	\norma{\dt p}_{ \L2 \Vp}
	+ \norma{\dt r}_{ \L2 \Vp}
	\leq c,
\end{align}
which, in turn, gives sense to the final conditions \eqref{Adj_IC}.
In fact, from the standard embedding of $\H1 \Vp \cap \L2 V $ in $\C0 H$,
we deduce that $p,r \in \C0 H.$

\noindent
{\bf Third estimate:} \,\,\,
Next, a comparison in \eqref{Adj_first} produces $\Delta q \in \L2 H$,
and the elliptic regularity theory yields that
\begin{align}
\label{adjoint_Third}
	\norma{q}_{ \L2 W}
	\leq c,
\end{align}
which also allows us to recover $q\in\C0 V$ from \wk~embedding results.

Summing up, we realize that the estimate
\begin{align}
\label{adjoint_final}
	\non & 
	\norma{q}_{ \H1 H \cap \C0 V \cap \L2 W}
	+ \norma{p}_{\H1 \Vp \cap \C0 H \cap \L2 V}
	\\[1mm]
	 & 
	+ \norma{r}_{\H1 \Vp \cap \C0 H \cap \L2 V}
	\leq c
\end{align}
has been proved.
The uniqueness part directly follows, since the
system \Adjsys\ is linear.
\Edim

Finally, we are left with the task of showing the necessary
conditions for optimality. To this end, we begin
by checking Theorem~\ref{THM_foc_prima}. Then, making use of
the adjoint system, we simplify \eqref{foc_first} and
deduce a variational inequality which is more convenient for the applications.

\Bdim[Proof of Theorem~\ref{THM_foc_prima}]
\juerg{This result is a direct consequence} of \eqref{formal_foc} and Theorem~\ref{THM_frechet}.
Indeed, combining the \Frechet\ differentiability of $\S$
with the chain rule, we can exploit \eqref{formal_foc} to derive \eqref{foc_first}.
\Edim

We are now in the position to eliminate the solutions to the linearized 
system from the necessary condition \eqref{foc_first}.
This procedure leads to \eqref{foc_final} and thus to Theorem~\ref{THM_foc_final}.

\begin{proof}[Proof of Theorem~\ref{THM_foc_final}]
Comparing the variational inequality
\eqref{foc_first} with \eqref{foc_final}, it becomes clear that
we only need to ensure that
\begin{align}
	& \non
	- \intQ h(\bph)q k 
	+\intQ lr
	=
	\bO \iO (\bph(T)-\phO)\xi(T)
	+ \bQ \intQ (\bph-\phQ)\xi
	\\ & \quad
	+ \bG \iO (\bs(T)-\sO)\z(T)
	+ \bS \intQ (\bs-\sQ)\z,
	\label{foc_proof}
\end{align}
where $\xi$ and $\z$ are the solution to the linearized
system \Linsys\ corresponding
to $k=u-\uopt$ and $l=w-\wopt$.
In order to show \eqref{foc_proof}, 
let us first point out that combining the Newton-Leibnitz formula
with the initial and final conditions \eqref{Lin_IC} and \eqref{Adj_IC}, respectively, we
have that
\begin{align*}
	\juerg{-\int_0^T \b \< \dt p(t), \xi(t)>_V\,dt}
	&\,=\, 
	\b \intQ \dt \xi \, p
	- \ioT\frac d {dt} \biggl( \iO \b p  \xi \biggr)\juerg{dt}\\
	& 
	=\,
	\b \intQ \dt \xi \, p
	- \iO \bO (\bph(T)-\phO)\xi(T),
	\\
	-\int_0^T  \juerg{\< \dt r(t), \z(t)>_V\,dt}
	&\,=\, 
	\intQ \dt \z \, r
	- \ioT\frac d {dt} \biggl( \iO r \z \biggr) \juerg{dt}
	\\
	&=\,
	\intQ \dt \z \, r
	- \iO \bG (\bs(T)-\sO) \z(T).
\end{align*}
Then, we consider the solution $(\et,\xi,\z)$ to \Linsys\ corresponding to
$k=u-\uopt$ and $l=w-\wopt$ as test functions in system \Adjsys. Namely,
we test \eqref{Adj_first} by $\et$, \eqref{Adj_second} by $\xi$,
\eqref{Adj_third} by $\z$, and integrate over $(0,T)$ to obtain that
\begin{align*} 
	0   &= \intQ \et\, [-\a \dt q - \Delta q- p]
		\\& \quad
		- \intQ \dt q\xi  	 	
		-\ioT \b \< \dt p(t),\xi(t)>_V\,dt  
		+  \intQ \nabla p \cdot \nabla \xi
		\andrea{-\ch \intQ \nabla r \cdot \nabla \xi}
		+ \intQ F''(\bph)p \xi
		\\ &\quad
		- \intQ (P\bs -A -\uopt)h'(\bph) q \xi
		+\intQ  D\bs h'(\bph) r \xi
		-\intQ \bQ (\bph-\phQ) \xi
	    \\ &\quad
	    -\ioT \<\dt r(t),\z(t)>_V\,dt + \intQ \nabla r \cdot \nabla\z 
	    +  \intQ Br \z
	    +\intQ D h(\bph)r\z
	    -\intQ \ch p \z
	   \\ &\quad
	    -\intQ P h(\bph)q \z
	    -\intQ \bS (\bs - \sQ) \z.
\end{align*}
\andrea{
Hence, we integrate by parts making use of the boundary conditions, the initial data and
the above identities. After rearrangements of the terms, we infer that
\begin{align*}
	&
	\iO \bO (\bph(T)-\phO)\xi(T)
	+\intQ \bQ (\bph-\phQ) \xi
	+ \iO \bG (\bs(T)-\sO) \z(T) 
	\\ & \quad 
	+\intQ \bS (\bs - \sQ)\z
	  =  \intQ p \,[	\b\dt \xi-\Delta \xi	+ F''(\bph)\xi	- \ch \z	-\et	]
		\\ & \quad
		+\intQ q \,[\a \dt \et	+ \dt \xi- \Delta \et - P\z h(\bph) -( P\bs - A -\uopt) h'(\bph)\xi]
	    \\ &\quad
	    + \intQ r\, [\dt \z	- \Delta \z	+B\z + \ch \Delta \xi+ D \z h(\bph)	+ D \bs h'(\bph)\xi].
\end{align*}
Finally,  we account for the equations of system \Linsys\ to realize that
\begin{align*}
	&
	\iO \bO (\bph(T)-\phO)\xi(T)
	+\intQ \bQ (\bph-\phQ) \xi
	+ \iO \bG (\bs(T)-\sO) \z(T) 
	\\ & \quad 
	+\intQ \bS (\bs - \sQ)\z
	  = 
		- \intQ h(\bph)q	(u-\uopt)
		+\intQ r(w-\wopt),
\end{align*}
that is \eqref{foc_proof}, so that the variational inequality \eqref{foc_final} has been shown.}
	 
Let us note, the last sentences in the statement of Theorem~\ref{THM_foc_final}
\sfw ly follow by combining the fact that
condition \eqref{foc_final} can be decoupled by taking 
first $w=\wopt$ and then $u=\uopt$ and use \juerg{the Hilbert projection theorem.}
\end{proof}

\section*{Acknowledgments}
\last{The authors are very grateful to the anonymous referee for the careful reading 
of the manuscript and for some useful suggestions.}
The research of Pierluigi Colli is supported by the Italian Ministry of Education,
University and Research~(MIUR): Dipartimenti di Eccellenza Program (2018--2022)
-- Dept.~of Mathematics ``F.~Casorati'', University of Pavia.
In addition, PC gratefully acknowledges some other
support from the GNAMPA (Gruppo Nazionale per l'Analisi Matematica,
la Probabilit\`a e le loro Applicazioni) of INdAM (Istituto
Nazionale di Alta Matematica) and the IMATI -- C.N.R. Pavia, Italy.

\vspace{3truemm}

\Begin{thebibliography}{10}
\footnotesize

%

%
\last{
\bibitem{BRZ}
H. Brezis,
``Op\'erateurs maximaux monotones et semi-groupes de contractions dans les
espaces de Hilbert'', \last{\it North-Holland Math. Stud.} {\bf 5}, North-Holland, Amsterdam, 1973.}

\bibitem{CRW}
C. Cavaterra, E. Rocca and H. Wu,
Long-time dynamics and optimal control of a diffuse interface model for tumor growth,
{\it Appl. Math. Optim.} (2019), \last{https://doi.org/10.1007/s00245-019-09562-5.}


\bibitem{CGH}
P. Colli, G. Gilardi and D. Hilhorst,
On a Cahn--Hilliard type phase field system related to tumor growth,
{\it Discrete Contin. Dyn. Syst.} {\bf 35} (2015), 2423-2442.

\bibitem{SM}
P. Colli, G. Gilardi, G. Marinoschi and E. Rocca,
Sliding mode control for a phase field system related to tumor growth,
{\it Appl. Math. Optim.} \textbf{79} (2019), 647-670.

%

\bibitem{CGRS_VAN}
P. Colli, G. Gilardi, E. Rocca and J. Sprekels,
Vanishing viscosities and error estimate for a Cahn--Hilliard type phase field system related to tumor growth,
{\it Nonlinear Anal. Real World Appl.} {\bf 26} (2015), 93-108.

\bibitem{CGRS_ASY}
P. Colli, G. Gilardi, E. Rocca and J. Sprekels,
Asymptotic analyses and error estimates for a Cahn--Hilliard type phase field system modelling tumor growth,
{\it Discrete Contin. Dyn. Syst. Ser. S} {\bf 10} (2017), 37-54.

\last{%
\bibitem{CGRS_OPT}
P. Colli, G. Gilardi, E. Rocca and J. Sprekels,
Optimal distributed control of a diffuse interface model of tumor growth,
{\it Nonlinearity} {\bf 30} (2017), 2518-2546.
\bibitem{CGS_frac_wp}
P. Colli, G. Gilardi and J. Sprekels, 
Well-posedness and regularity for a fractional tumor growth model,
{\it Adv. Math. Sci. Appl.} {\bf 28} (2019), 343-375.
\bibitem{CGS_frac_control}
P. Colli, G. Gilardi and J. Sprekels, 
A distributed control problem for a fractional tumor growth model,
{\it Mathematics} {\bf 7} (2019), https://doi:10.3390/math7090792.
\bibitem{CGS_frac_asy}
P. Colli, G. Gilardi and J. Sprekels, 
Asymptotic analysis of a tumor growth model with fractional operators,
preprint arXiv:1908.10651 [math.AP] (2019), 1-31.}

\bibitem{CLLW}
V. Cristini, X. Li, J.S. Lowengrub and S.M. Wise,
Nonlinear simulations of solid tumor growth using a mixture model: invasion and branching.
{\it J. Math. Biol.} {\bf 58} (2009), 723-763.

\bibitem{CL}
V. Cristini and J. Lowengrub,
``Multiscale Modeling of Cancer: An Integrated Experimental and Mathematical
Modeling Approach'', Cambridge University Press, Leiden (2010).

\bibitem{DFRGM}
M. Dai, E. Feireisl, E. Rocca, G. Schimperna and \last{M.E.} Schonbek,
Analysis of a diffuse interface model of multi-species tumor growth,
{\it Nonlinearity\/} {\bf  30} (2017), 1639-1658.

\bibitem{EK_ADV}
M. Ebenbeck and P. Knopf,
Optimal control theory and advanced optimality conditions 
for a diffuse interface model of tumor growth,
preprint arXiv:1903.00333 [math.OC] (2019), 1-34.

\last{
\bibitem{EK}
M. Ebenbeck and P. Knopf,
Optimal medication for tumors modeled by a Cahn--Hilliard--Brinkman equation,
\last{{\it Calc. Var. Partial Differential Equations}} {\bf 58} (2019), 
https://doi.org/10.1007/s00526-019-1579-z.%
}

\bibitem{EGAR}
M. Ebenbeck and H. Garcke,
Analysis of a Cahn--Hilliard--Brinkman model for tumour growth with chemotaxis,
{\it J. Differential Equations} \textbf{266} (2019), 5998-6036.

\bibitem{FGR}
S. Frigeri, M. Grasselli and E. Rocca,
On a diffuse interface model of tumor growth,
{\it  European J. Appl. Math.\/} {\bf 26} (2015), 215-243. 

\bibitem{FLRS}
S. Frigeri, K.F. Lam, E. Rocca and G. Schimperna,
On a multi-species Cahn--Hilliard--Darcy tumor growth model with singular potentials,
{\it Commun. Math Sci.} \last{{\bf  16}} (2018), 821-856. 

\bibitem{FRL}
S. Frigeri, K.F. Lam and E. Rocca,
On a diffuse interface model for tumour growth with non-local 
interactions and degenerate mobilities, 
in ``Solvability, regularity, and optimal control of boundary 
value problems for PDEs'', 
P.~Colli, A.~Favini, E.~Rocca, G.~Schimperna, J.~Sprekels~(ed.),
{\it Springer INdAM Series}~{\bf 22}, Springer, Cham, 2017, 217-254.

%
\bibitem{GARL_3}
H. Garcke and K.F. Lam,
Global weak solutions and asymptotic limits 
of a Cahn--Hilliard--Darcy system modelling tumour growth,
{\it AIMS Mathematics} {\bf 1} (2016), 318-360.
\bibitem{GARL_1}
H. Garcke and K.F. Lam,
Well-posedness of a Cahn--Hilliard system modelling tumour
growth with chemotaxis and active transport,
{\it European. J. Appl. Math.} {\bf 28} (2017), 284-316.
\bibitem{GARL_2}
H. Garcke and K.F. Lam,
Analysis of a Cahn--Hilliard system with non--zero Dirichlet 
conditions modeling tumor growth with chemotaxis,
{\it Discrete Contin. Dyn. Syst.} {\bf 37} (2017), 4277-4308.
\bibitem{GARL_4}
H. Garcke and K.F. Lam,
On a Cahn--Hilliard--Darcy system for tumour growth 
with solution dependent source terms, 
in ``Trends on Applications of Mathematics to Mechanics'', 
E.~Rocca, U.~Stefanelli, L.~Truskinovski, A.~Visintin~(ed.), 
{\it Springer INdAM Series} {\bf 27}, Springer, Cham, 2018, 243-264.
\bibitem{GAR}
H. Garcke, K.F. Lam, R. N\"urnberg and E. Sitka,
A multiphase Cahn--Hilliard--Darcy model for tumour growth with necrosis,
{\it Math. Models Methods Appl. Sci.} {\bf 28} (2018), 525-577.

\bibitem{GARLR}
H. Garcke, K.F. Lam and E. Rocca,
Optimal control of treatment time in a diffuse interface model of tumor growth,
{\it Appl. Math. Optim.} {\bf 78} (2018), {495-544}.

\bibitem{GLSS}
H. Garcke, K.F. Lam, E. Sitka and V. Styles,
A Cahn--Hilliard--Darcy model for tumour growth with chemotaxis and active transport,
{\it Math. Models Methods Appl. Sci. } {\bf 26} (2016), 1095-1148.


\bibitem{HDPZO}
A. Hawkins-Daarud, S. Prudhomme, K.G. van der Zee \last{and} J.T. Oden,
Bayesian calibration, validation, and uncertainty quantification of diffuse 
interface models of tumor growth. 
{\it J. Math. Biol.} {\bf 67} (2013), 1457-1485. 

\bibitem{HDZO}
A. Hawkins-Daarud, K.G. van der Zee and J.T. Oden, Numerical simulation of
a thermodynamically consistent four-species tumor growth model, 
{\it Int. J. Numer. Math. Biomed. Engng.} {\bf 28} (2011), 3-24.

\bibitem{HKNZ}
D. Hilhorst, J. Kampmann, T.N. Nguyen and K.G. van der Zee, 
Formal asymptotic limit of a diffuse-interface tumor-growth model, 
{\it Math. Models Methods Appl. Sci.} {\bf 25} (2015), 1011-1043.


\bibitem{LDY}
O.A. Lady\v zenskaja, V.A. Solonnikov and N.N. Uralceva,
``Linear and quasilinear equations of parabolic type'',
{\it Mathematical Monographs Volume} {\bf 23},
American Mathematical Society, Providence, 1968.

\bibitem{Lions_OPT}
J.-L. Lions,
``Contr\^ole optimal de syst\`emes gouvern\'es par des equations aux d\'eriv\'ees partielles'',
Dunod, Paris, 1968.




\bibitem{MRS}
A. Miranville, E. Rocca and G. Schimperna,
On the long time behavior of a tumor growth model,
{\it J. Differential Equations\/} {\bf 267} (2019)\last{,} 2616-2642.
%
%

\bibitem{OHP}
J.T. Oden, A. Hawkins and S. Prudhomme,
General diffuse-interface theories and an approach to predictive tumor growth modeling,
{\it Math. Models Methods Appl. Sci.} {\bf 20} (2010), 477-517. 

\last{%
\bibitem{S}
A. Signori,
Optimal distributed control of an extended model of tumor 
growth with logarithmic potential,
{\it Appl. Math. Optim.} \last{(2019)}, https://doi.org/10.1007/s00245-018-9538-1.
\bibitem{S_DQ}
A. Signori,
Optimality conditions for an extended tumor growth model with 
double obstacle potential via deep quench approach, 
{\it Evol. Equ. Control Theory} (2019), https://doi: 10.3934/eect.2020003.
\bibitem{S_b}
A. Signori,
Optimal treatment for a phase field system of Cahn--Hilliard 
type modeling tumor growth by asymptotic scheme, 
{\it Math. Control Relat. Fields} (2019), https://doi: 10.3934/mcrf.2019040.
\bibitem{S_a}
A. Signori,
Vanishing parameter for an optimal control problem modeling tumor growth,
preprint arXiv:1903.04930 [math.AP] (2019), 1-22.}

\bibitem{Simon}
J. Simon,
{Compact sets in the space $L^p(0,T; B)$},
{\it Ann. Mat. Pura Appl.~(4)\/} 
{\bf 146} (1987), 65-96.

\bibitem{SW}
J. Sprekels and H. Wu,
Optimal distributed control of a Cahn--Hilliard--Darcy system with mass sources,
{\it Appl. Math. Optim.} (2019), https://doi.org/10.1007/s00245-019-09555-4.

\bibitem{Trol}
F. Tr\"oltzsch,
``Optimal Control of Partial Differential Equations. Theory, Methods and Applications'',
{\it Grad. Stud. in Math.} {\bf 112}, AMS, Providence, RI, 2010.

\bibitem{WLFC}
S.M. Wise, J.S. Lowengrub, H.B. Frieboes and V. Cristini,
Three-dimensional multispecies nonlinear tumor growth - I: model and numerical method, 
{\it J. Theor. Biol.} \last{{\bf 253}} (2008), 524-543.

\bibitem{WZZ}
X. Wu, G.J. van Zwieten and K.G. van der Zee, Stabilized second-order splitting
schemes for \CH~models with applications to diffuse-interface tumor-growth models, 
{\it Int. J. Numer. Meth. Biomed. Engng.} {\bf 30} (2014), 180-203.

\End{thebibliography}

\End{document}

\bye